\documentclass[12pt,reqno]{amsart}
\usepackage{amsmath,amsfonts,amsthm,color,amssymb, mathrsfs}

\usepackage{pxfonts}

\makeatletter
\@tfor\@tempa:=
\alpha \beta \gamma \delta \epsilon \zeta \eta \theta \iota \kappa \lambda \mu
\nu \xi \pi \rho \sigma \tau \upsilon \phi \chi \psi \omega \varepsilon
\vartheta \varpi \varrho \varsigma \varphi  
\do {
	\count255 \numexpr\@tempa+"7000\relax
	\expandafter\mathchardef\@tempa \count255 }
\makeatother

\usepackage{tikz-cd}
\usetikzlibrary{arrows.meta,decorations.pathreplacing,calc}
\usepackage{xcolor}
\usepackage[T1]{fontenc}
\usepackage{wasysym}
\usepackage[normalem]{ulem}
\usepackage{stmaryrd} 

\usepackage[left=1in, right=1in, top=1.1in,bottom=1.1in]{geometry}
\usepackage{enumitem}
\usepackage{hyperref}
\hypersetup{
	colorlinks   = true,
	citecolor    = blue,
	linkcolor=blue
}
\allowdisplaybreaks
\usepackage{csquotes}
\usepackage{graphicx}
\usepackage{pstricks}
\usepackage{lmodern}

\newtheorem{lemma}{Lemma}
\newtheorem{theorem}{Theorem}[section]

\newtheorem{definition}{Definition}

\newtheorem{proposition}{Proposition}

\newtheorem{remark}{Remark}

\newcommand{\bean}{\begin{eqnarray*}}
	\newcommand{\eean}{\end{eqnarray*}}

\newcounter{bean}
\newcommand{\benuma}{\setlength{\labelwidth}{.25in}
	
	\begin{list}
		{(\alph{bean})}{\usecounter{bean}}}
	\newcommand{\eenuma}{\end{list}}

\begin{document}
	
	\title[The most probable transition paths in multiscale systems]{Convergence of the Most Probable Transition Paths in Multiscale Systems}
	
	\author[Q. Lian]{Qiqi Lian}
	\address{School of Mathematics, Sun Yat-sen University, Guangzhou, Guangdong, 510275, China}
	\address{School of Sciences, Great Bay University, Dongguan, Guangdong, 523000, China}
	\email{lianqq3@mail2.sysu.edu.cn}
	
	\author[S. Liu]{Shanqi Liu}
	\address{School of Sciences, Great Bay University, Dongguan, Guangdong, 523000, China}
	\address{School of Mathematical Sciences, University of Science and Technology of China, Hefei, Anhui, 230026, China}
	\email{shanqiliumath@126.com}
	
	\author[J. Duan]{Jinqiao Duan}
	\address{School of Sciences, Great Bay University, Dongguan, Guangdong, 523000, China}
	\email{duan@gbu.edu.cn (Corresponding author), ORCID: 0000-0002-2077-990X}
	
	\vspace{-2cm}
	\maketitle
	\vspace{-1cm}
	\begin{abstract}		
  	This paper is devoted to study the convergence behavior of the most probable transition paths in a class of multiscale stochastic dynamical systems. 
	First, based on Euler-Lagrange equations, we establish proximity and convergence results in the space of continuous functions endowed with the supremum norm.
	These results compare the most probable paths of the multiscale system with those of averaged system derived from the averaging principle.
	Next, we utilize a variational framework and analyze the Onsager-Machlup functionals using $\Gamma$-convergence.
	Together with the equi-coercivity of the functional family, we prove the convergence of both their minima and minimizers.
	Specifically, we show weak convergence of the minimizers in the Cameron-Martin space $H^1$, which implies their strong convergence in $L^2$.
		
		
	\end{abstract}
	\section{Introduction}
	Stochastic dynamical systems provide a fundamental framework for modeling a wide range of complex phenomena influenced by random fluctuations. 
	In many such systems, particularly those with multiple stable states, noise plays a crucial role in transitions between metastable states. 
	Such transitions are significant to understanding diverse processes, including climate tipping points, protein folding pathways, chemical reaction kinetics, and the onset of epileptic seizures \cite{hanggi1990reaction, lenton2008tipping, scheffer2009early, takeshita2007transitions}.
	One important issue is to determine the \textit{most probable transition path} (MPTP), which describes the most likely trajectory of the system and 
	establishes a connection between deterministic dynamics and stochastic fluctuations.
 
	A typical approach to characterizing the MPTP is through the Onsager-Machlup (OM) action functional via a variational framework. 
	Originally introduced in 1953 within the context of non-equilibrium statistical mechanics by Onsager and Machlup \cite{onsager1953fluctuations, machlup1953fluctuations}, the OM functional initially served as the probability density functional for diffusion processes characterized by linear drift and constant diffusion coefficients. 
	Subsequently, Tisza and Manning \cite{tisza1957fluctuations} generalized these foundational results to nonlinear equations.	
	Stratonovich \cite{stratonovich1971probability} provided a rigorous mathematical justification. 
	Since then, OM functionals have been derived for a variety of stochastic processes or measures \cite{capitaine2007onsager, carfagnini2024onsager, chao2019onsager,
		coulibaly2014onsager, durr1978onsager, fujita1982onsager, gasteratos2026onsager, huang2025probability, ikeda2014stochastic, liu2023onsager, moret2002onsager, takahashi2006probability}.
	In general, let $X_t$ be an $L_t$-diffusion process starting at point $x_0$ on a Riemannian manifold $N$, where the generator is defined as 
	$L_t = \frac{1}{2}\Delta_t + Z(t, \cdot)$. 
	Then, for a smooth curve (reference path) $\xi(t)$, 
	the following asymptotic behavior of the tube holds \cite{coulibaly2014onsager}
	\begin{equation*}
		\mathbb{P}_{x_0}\left\{\forall t \in [0,T], \, d(t, X_t, \xi(t)) \leq \delta \right\} 
		\sim_{\delta \downarrow 0}
		C \exp\left\{-\frac{\lambda_1 T}{\delta^2}\right\} 
		\exp\left\{-\frac{1}{2}\int_0^T M(t, \xi(t), \dot{\xi}(t))  \mathrm{d}t \right\},
	\end{equation*}
	where $d(t, x, y)$ is the Riemannian distance on $N$ according to the metric $g(t)$, 
	$M$ is a time-dependent function on the tangent bundle, known as the \textit{OM function} and defined for $v \in T_xN$ as 
	\begin{equation*}
		M(t, x, v) = \|Z(t, x) - v\|_{g(t)}^2 +\operatorname{div}_{g(t)}(Z)(t, x) - \frac{1}{6} R_{g(t)}(x) + \frac{1}{2} \operatorname{trace}_{g(t)}(\dot{g}(t)).
	\end{equation*}
	Here, $C$ and $\lambda_1$ are explicit constants, while $\operatorname{div}_{g(t)}$ and $R_{g(t)}$ represent the divergence operator and the scalar curvature with respect to the time-dependent metric $g(t)$, respectively. 
	The integral $\int_0^T M(t, \varphi(t), \dot{\varphi}(t)) \mathrm{d}t$ is known as the \textit{OM functional}.
    The OM function can be regarded as a Lagrangian, with the corresponding OM functional serving as an action functional \cite{durr1978onsager, graham1977path, tisza1957fluctuations}. 
     Hence, the most probable transition path can be characterized by the Euler-Lagrange equation.

However, applying this variational framework directly to multiscale systems is considerably complicated by their inherent separation of scales. 
The resulting Euler-Lagrange equations inevitably inherit the multiscale structure of the original dynamics, typically exhibiting severe nonlinearities and rapid oscillations, 
which makes their asymptotic analysis extremely challenging. 
To circumvent the prohibitive cost of directly simulating the full multiscale system, the averaging principle is widely adopted as an effective model reduction technique. 
	Following the seminal works of Krylov, Bogoliubov, Khasminskii, Freidlin, and others 
	(see \cite{bogoliubov1961asymptotic, cerrai2009averaging, duan2014effective, freidlin1998random, gao2022averaging, khas1968averaging, 
		krylov2016introduction, pavliotis2008multiscale, 
		weinan2007heterogeneous} and references therein), 
	the averaging principle allows one to systematically eliminate fast-varying components. Consequently, one can derive a simplified, effective stochastic differential equation (SDE) that exclusively governs the dynamics of the slow variables. 
	This averaging principle naturally raises an interesting question: can the most probable transition path of the averaged SDE serve as a reliable substitute for the most probable transition path of the original multiscale system? 
	Specifically, it is crucial to establish whether a rigorous approximation or convergence relationship exists between the transition paths as the scale parameter tends to zero.

    In addition to the Euler-Lagrange characterization of the most probable transition paths, we can derive convergence relations between the minimizers by directly studying the convergence of the Onsager-Machlup functionals.
	An effective tool for this purpose is \textit{$\Gamma$-convergence}, 
	which is a commonly-recognized notion of convergence for variational problems.
	The $\Gamma$-convergence of OM functionals has already been studied in some cases. 
	For example, Pinski et al. \cite{pinski2012gamma} discussed the $\Gamma$-convergence of the OM functional at low temperatures, which corresponds to the small-noise regime. 
	Similarly with the small noise, Li and Li \cite{li2021gamma} demonstrated that the $\Gamma$-limit of the OM functional on the space of curves reduces to the geometric form of the Freidlin-Wentzell functional.
	In our setting, however, we consider finite noise rather than the small-noise regime. 
	
	Inspired by all these above, we aim to investigate the relationship between the most probable transition path in multiscale system and that derived from the effective system under the Bogoliubov averaging principle.
	Our main results are two folds and are illustrated in the diagram below.
	First, based on Euler-Lagrange equations, we establish the proximity result (Theorem \ref{main}) and the convergence result (Theorem \ref{mainn}) of the transition paths in the space of continuous functions endowed with the supremum norm. 
	Furthermore, by establishing $\Gamma$-convergence and the equi-coerciveness of the family of OM functionals, we prove the weak convergence of minima and minimizers (Theorem \ref{thm:convergence_minima}) in the Cameron-Martin space $H^1$ and their consequent strong convergence in the $L^2$ topology. These results provide a theoretical foundation for utilizing reduced averaged models to analyze rare events.
	
	\begin{figure}[htbp]
		\centering
		\includegraphics[width=0.9\textwidth]{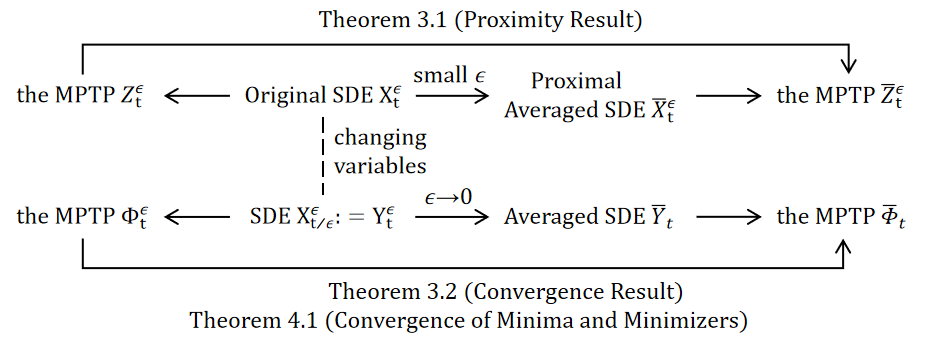}
	\end{figure}
	
	The rest of this paper is organized as follows. In Section \ref{Pre}, we introduce preliminaries and assumptions about the averaging principle, the OM action functional, the most probable transition path, and $\Gamma$-convergence. We prove the proximity and convergence of the most probable transition paths based on Euler-Lagrange equations in Section \ref{PC}. In Section \ref{GC}, we utilize $\Gamma$-convergence for the OM functionals to rigorously prove the convergence of minima and their respective minimizers. In Section \ref{NE}, we present a concrete example to illustrate our results. 
	Section \ref{DS} is a discussion that summarizes the contributions of our work.
	\section{Preliminaries}\label{Pre}
	In this section, we present the framework and preliminaries needed throughout the paper. We begin with the averaging principle and the assumptions on the coefficients. 
	Next, we recall the Onsager-Machlup action functional, the most probable transition path, and the theory of $\Gamma$-convergence. 
	We conclude this section with a generalized Gronwall-Bellman-Bihari inequality that will play a key role in the subsequent analysis.
	
	\subsection{Averaging principle}
	Let $(\Omega, \mathcal{F}, \mathbb{P})$ be a completed probability space. 
	Consider the following stochastic differential equation in $\mathbb{R}^d$
	\begin{align}\label{SDE1}
		\mathrm{d}X^{\epsilon}_t=\epsilon b(t,X^{\epsilon}_t)\mathrm{d}t+\sqrt{\epsilon} \mathrm{d}B_t,
	\end{align}
	where $t\in[0,1]$, the initial value is $X^{\epsilon}_0=x_0\in\mathbb{R}^d$, the terminal value is $X^{\epsilon}_1=x_1\in\mathbb{R}^d$, 
	$0<\epsilon<1$ is a small parameter, the coefficient $b:[0,1]\times \mathbb{R}^d\rightarrow \mathbb{R}^d$ is a measurable function
	and $B_t$ is a standard $d$-dimensional Brownian motion. 
	Furthermore, by changing variables, we can obtain 
	\begin{align}
		X^{\epsilon}_{t/\epsilon}
		&=x_0+\epsilon\int_0^{t/\epsilon}  b(s,X^{\epsilon}_s)\mathrm{d}s
		+\sqrt{\epsilon}\int_0^{t/\epsilon}  \mathrm{d}B_s,\nonumber\\
		&=x_0+\int_0^t b({s/\epsilon},X^{\epsilon}_{s/\epsilon})\mathrm{d}s+ B_t.\nonumber 
	\end{align}
	Denote $X^{\epsilon}_{t/\epsilon}:=Y_t^\epsilon$, it follows that
	\begin{align*}
		Y_t^\epsilon
		=x_0+\int_0^t b({s/\epsilon},Y_s^\epsilon)\mathrm{d}s+ B_t,
	\end{align*}
	that is 
	\begin{align} \label{SDE11} 
		\mathrm{d}Y_t^\epsilon
		=b({t/\epsilon},Y_t^\epsilon)\mathrm{d}t+ \mathrm{d}B_t,
	\end{align}
	which is equivalent to Equ.\eqref{SDE1} after changing variables.
		
	The averaged SDEs for the original SDE \eqref{SDE1} are 
	\begin{align}\label{SDE2}
		\mathrm{d}\bar{X}^{\epsilon}_t=\epsilon \bar{b}(\bar{X}^{\epsilon}_t)\mathrm{d}t+\sqrt{\epsilon}\mathrm{d}B_t,
	\end{align}
	and
	\begin{align}\label{SDE3}
		\mathrm{d}\bar{Y}_t=\bar{b}(\bar{Y}_t)\mathrm{d}t+\mathrm{d}B_t,
	\end{align}
	where the initial value $\bar{X}^{\epsilon}_0=\bar{Y}_0=x_0\in\mathbb{R}^d$  
	and the coefficient $\bar{b}:\mathbb{R}^d\rightarrow \mathbb{R}^d$ is a measurable function. 
	We need to impose the following assumptions on the coefficients.
	\begin{itemize}
		\item[$\bullet$ \textbf{(A1)}] 
		Lipschitz condition. There exists a bounded function $K_1(t)>0$ such that for all $y_1, y_2\in\mathbb{R}^d$,
		\begin{align}\label{b1}
			|b(t,y_1)-b(t,y_2)|\leq K_1(t)|y_1-y_2|.
		\end{align}
		Growth condition. There exists a bounded function $K_2(t)>0$ such that for all $y\in\mathbb{R}^d$,
		\begin{align}\label{b2}
			|b(t,y)|\leq K_2(t)(1+|y|).
		\end{align}
	\end{itemize}

	\begin{itemize}
		\item[$\bullet$ \textbf{(A2)}] 
		Averaged condition. 	
		There exists a bounded function $\phi_1(T)>0$ with $\lim\limits_{T_1\rightarrow \infty}\phi_1(T_1)=0$, 
		such that for all $y\in\mathbb{R}^d$ and $T\in[0,1]$,
		\begin{align}\label{ac}
			\frac{1}{T}\int_{0}^{T}	|b(t,y)-\bar{b}(y)|\mathrm{d} t\leq \phi_1(T)(1+|y|).
		\end{align}
	\end{itemize}
	\begin{remark}
		(i) From the averaged condition, it can be deduced that 
		$$\lim\limits_{T_1\rightarrow \infty}\frac{1}{T_1}
		\int_{0}^{T_1}b(t,y)\mathrm{d} t=\bar{b}(y),$$
		where the convergence is assumed uniform with respect to $y$. 
		
		(ii) Together with the assumptions $\textbf{(A1)}$ and $\textbf{(A2)}$, 
		we can deduce the Lipschitz and growth properties of the averaged coefficient $\bar{b}$ also hold, i.e. 
			there exist positive constants $K_3$ and $K_4$ such that for all $y_1, y_2\in\mathbb{R}^d$,
			\begin{align}\label{bb1}
				|\bar{b}(y_1)-\bar{b}(y_2)|\leq K_3|y_1-y_2|, 
			\end{align}
			and
			\begin{align}\label{bb2}
				|\bar{b}(y_1)|\leq K_4(1+|y_1|).
			\end{align}
	\end{remark}
	
	According to \cite{freidlin1998random} or \cite{xu2011averaging}, 
	we know that the solution to Equ.\eqref{SDE2} approximates the solution to Equ.\eqref{SDE1} when $\epsilon$ is sufficiently small, and the solution to Equ.\eqref{SDE11} converges in probability to the solution to Equ.\eqref{SDE3} as $\epsilon$ tends to zero. 
	More precisely, the following two lemmas hold.
	\begin{lemma}[\cite{xu2011averaging}]\label{L2}
		Suppose that the assumptions \textbf{(A1)–(A2)} are satisfied. For a given arbitrarily small number $\delta_1 > 0$, there exist $L > 0$, $\epsilon_1 \in (0, \epsilon_0)$ and $\beta \in (0, 1)$ such that for all $\epsilon \in (0, \epsilon_1]$, 
		$$\mathbb{E} \left(\sup_{t \in [0, L\epsilon^{-\beta}]} |X^\epsilon_t - \bar{X}^\epsilon_t|^2 \right)\leq \delta_1.$$ 
	\end{lemma}
	\begin{lemma}[\cite{xu2011averaging}]\label{P}
		Suppose that the assumptions \textbf{(A1)–(A2)} are satisfied. Then for any number $\delta_2 > 0$ and any finite $t_0$, we have 
		$$\lim_{\epsilon \to 0} \mathbb{P} \left( \sup_{t \in (0, t_0]} |X^\epsilon_{t/\epsilon} - \bar{Y}_t| < \delta_2 \right) = 1.$$ 
	\end{lemma}
	\subsection{Onsager-Machlup action functional}
	
	Let the time interval be $[0,1]$. The space of paths for a diffusion process starting at $x_0$ is denoted by $C_{1}^{x_{0}}$, which is the space of continuous functions mapping $[0,1]$ to $\mathbb{R}^d$:
	\begin{equation*}
		C_{1}^{x_{0}} = \left\{ x : [0,1] \rightarrow \mathbb{R}^d \mid x(t) \text{ is continuous, } x(0) = x_0 \right\}.
	\end{equation*}
	Equipped with the uniform norm $\|\cdot\|$, defined as
	\begin{equation*}
		\|x\| = \max_{t \in [0,1]} |x(t)|, \quad \text{for } x \in C_{1}^{x_{0}},
	\end{equation*}
	the space $C_{1}^{x_{0}}$ is a Banach space. 
	
	We require that the shifted path $\xi - x_0$ belongs to the Cameron-Martin space $H^1$. Here, $H^1$ denotes the class of all absolutely continuous functions $h$ such that $h(0) = 0$ and $\dot{h} \in L^2([0, 1]; \mathbb{R}^d)$. The inner product on $H^1$ is given by the formula
	\begin{equation*}
		(h_1, h_2)_{H^1} \coloneqq \int_0^1 \dot{h}_1(t)\dot{h}_2(t)\,\mathrm{d}t.
	\end{equation*}
	
	\begin{definition}[Onsager-Machlup Functional]
		Let $\delta > 0$ be given. Consider a diffusion process $X$ on $[0,1]$ and 
		a reference path $\xi$. 
		For a sufficiently small $\delta$, if the probability of the process $X$ lying in the tube can be expressed in the asymptotic form, 
		$$
		\mathbb{P}\left(\|X - \xi\| \leq \delta\right)
		\propto C(\delta) \exp\left\{-\frac{1}{2} \int_0^1 M(\xi,\dot{\xi})\,\mathrm{d}t\right\},
		$$
		then the integrand $M(\xi,\dot{\xi})$ is called the \textit{Onsager-Machlup function}. Here, $\propto$ denotes asymptotic equivalence as $\delta \to 0$. We refer to the integral $\int_0^1 M(\xi,\dot{\xi})\,\mathrm{d}t$ as the \textit{Onsager-Machlup functional}. 
		Furthermore, OM functionals are unique only up to the addition of real constants. 
	\end{definition}
	
	According to classical results (e.g.\cite{durr1978onsager, ikeda2014stochastic}), 
	we can explicitly derive the OM functions for the four distinct diffusion processes considered in this study. Let $z^\epsilon$, $\varphi^\epsilon$, $\bar{z}^\epsilon$ and $\bar{\varphi}$ denote the respective reference paths for the solutions to Eqs.\eqref{SDE1}, \eqref{SDE11}, \eqref{SDE2} and \eqref{SDE3}. 
	We present their corresponding OM functions collectively:
	\begin{align}\label{OM1}
		M_1(\dot{z}^\epsilon,z^\epsilon)
		=\left(\frac{\epsilon b(t,z^{\epsilon}_t)-\dot{z}^\epsilon_t}{\sqrt{\epsilon}}\right)^2
		+\epsilon b'(t,z^{\epsilon}_t),
	\end{align}
	\begin{align}\label{OM11}
		M_{2}(\dot{\varphi}^\epsilon,\varphi^\epsilon)
		=\left(b(t/\epsilon,\varphi^{\epsilon}_t)-\dot{\varphi}^\epsilon_t\right)^2
		+b'(t/\epsilon,\varphi^{\epsilon}_t),
	\end{align}
	\begin{align}\label{OM2}
		M_3(\dot{\bar{z}}^\epsilon,\bar{z}^\epsilon)
		=\left(\frac{\epsilon \bar{b}(\bar{z}^{\epsilon}_t)-\dot{\bar{z}}^\epsilon_t}{\sqrt{\epsilon}}\right)^2
		+\epsilon\bar{b}'(\bar{z}^{\epsilon}_t),
	\end{align}
	and
	\begin{align}\label{OM22}
		M_{4}(\dot{\bar{\varphi}},\bar{\varphi})
		=\left(\bar{b}(\bar{\varphi}_t)-\dot{\bar{\varphi}}_t\right)^2
		+\bar{b}'(\bar{\varphi}_t),
	\end{align}
	where $b'(t,y):=\partial_x b(t,x)\big|_{x=y}$ 
	denotes the partial derivatives with respect to the spatial variable $x$.
	\subsection{The most probable transition path}
	The minimization of the Onsager-Machlup functional naturally leads to the notion of the most probable transition path.
	
	\begin{definition}[The Most Probable Transition Path]\label{def:mptp}
		The \textit{most probable transition path} $\Xi$ is defined as the minimizer of the OM action functional $\int_{0}^{1}M(\xi, \dot{\xi})\mathrm{d}t$. 
		It is the solution to the following variational problem, subject to pin-down conditions $x_0$ and $x_1$ at the boundaries of the time interval $[0,1]$
		\begin{equation*}  
			\Xi 
			=\operatorname*{arg\,min}_{\substack{\xi \in H^1(0,1) \\ \xi_0 = x_0, \, \xi_1 = x_1}}
			\int_0^1 M(\xi, \dot{\xi}) \mathrm{d}t.
		\end{equation*}
	\end{definition}
   From a classical mechanics perspective, the calculus of variations together with Tonelli's theorem \cite{fathi2008weak} shows that
	the above problem has a solution $\Xi$, which are restricted in $C^2([0,1];\mathbb{R}^d)$ and satisfies the Euler-Lagrange equation 
	\begin{equation*} 
		\delta \int_0^1 M(\xi, \dot{\xi})\mathrm{d}t = 0,
	\end{equation*}
	that is,
	\begin{equation*} 
		\frac{\mathrm{d}}{\mathrm{d}t} 
		\frac{\partial M(\Xi, \dot{\Xi})}
		{\partial \dot{\Xi}} 
		= \frac{\partial M(\Xi, \dot{\Xi})}{\partial \Xi}.
	\end{equation*}
	Under certain conditions for an SDE with drift $b(t, x)$ and diffusion $\sigma(t,x)=1$, this leads to the following specific second-order ODE
	\begin{align*} 
		\left\{
		\begin{array}{ll}
			\ddot{\Xi}_t = b(t,\Xi_t) b'(t,\Xi_t) +\frac{1}{2} b''(t,\Xi_t)
			+ \partial_t b(t,\Xi_t), \\[6pt]
			\Xi_0=  x_0, \quad \Xi_1 = x_1,
		\end{array}
		\right.
	\end{align*}
	where $b'(t,\Xi_t):=\partial_x b(t,x)\big |_{x = \Xi_t}$, $b''(t,\Xi_t):=\partial_{x}^2 b(t,x)\big |_{x = \Xi_t}$, and $\partial_t b(t,x)$ denote the partial derivatives with respect to $x$ and $t$, respectively.
	
	Consequently, given the OM functions defined in \eqref{OM1}, \eqref{OM11}, \eqref{OM2} and \eqref{OM22}, it follows that the corresponding most probable transition paths $Z^{\epsilon}$, $\varPhi^{\epsilon}$, $\bar{Z}^{\epsilon}$ and $\bar{\varPhi}$ satisfy the following Euler-Lagrange equations:
	\begin{equation}\label{EL1}
		\ddot{Z}^\epsilon_t = \epsilon^2 b(t,Z^{\epsilon}_t)b'(t,Z^{\epsilon}_t)
		+\frac{  \epsilon^2}{2} b''(t,Z^{\epsilon}_t)
		+\epsilon\partial_t b(t,Z^{\epsilon}_t), 
	\end{equation}
	\begin{equation}\label{EL11}
		\ddot{\varPhi}^\epsilon_t = b(t/\epsilon,\varPhi^{\epsilon}_t)b'(t/\epsilon,\varPhi^{\epsilon}_t)
		+\frac{1}{2} b''(t/\epsilon,\varPhi^{\epsilon}_t)
		+\frac{1}{\epsilon}\partial_t b(t/\epsilon,\varPhi^{\epsilon}_t), 
	\end{equation}
	\begin{equation}\label{EL2}
		\ddot{\bar{Z}}^\epsilon_t = \epsilon^2 \bar{b}(\bar{Z}^{\epsilon}_t)\bar{b}'(\bar{Z}^{\epsilon}_t)
		+\frac{  \epsilon^2}{2} \bar{b}''(\bar{Z}^{\epsilon}_t), 
	\end{equation}
	and
	\begin{equation}\label{EL3}
		\ddot{\bar{\varPhi}}_t = 
		\bar{b}(\bar{\varPhi}_t)\bar{b}'(\bar{\varPhi}_t)
		+\frac{1}{2} \bar{b}''(\bar{\varPhi}_t),
	\end{equation}
	with boundary conditions
	\begin{equation}
		Z^\epsilon_0=\varPhi^\epsilon_0=\bar{Z}^{\epsilon}_0=\bar{\varPhi}_0= x_0, \quad 
		Z^\epsilon_1=\varPhi^\epsilon_1=\bar{Z}^{\epsilon}_1=\bar{\varPhi}_1= x_1. 
	\end{equation}
	Since the derived Euler-Lagrange equations involve the spatial derivatives $b'$, $b''$, $\bar{b}'$ and $\bar{b}''$, as well as the temporal derivative $\partial_t b$, further regularity assumptions on the drift coefficients are required.
	\begin{itemize}
		\item[$\bullet$ \textbf{(A3)}] 
		Lipschitz condition. There exists a bounded function $K_5(t)>0$ such that for all $y_1, y_2\in\mathbb{R}^d$,
		\begin{align}\label{b'1}
			|b'(t,y_1)-b'(t,y_2)|+|b''(t,y_1)-b''(t,y_2)|\leq K_5(t)|y_1-y_2|.
		\end{align}	
		Growth condition. There exists a bounded function $K_6(t)>0$ such that for all $y\in\mathbb{R}^d$,
		\begin{align}\label{b'2}
			|\partial_t b(t,y)|+|b'(t,y)|+|b''(t,y)|\leq K_6(t)(1+|y|).
		\end{align}
		\item[$\bullet$ \textbf{(A4)}] 
		Averaged condition. 
		There exists a bounded function $\phi_2(T)>0$ with $\lim\limits_{T_1\rightarrow \infty}\phi_2(T_1)=0$, such that for all $y\in\mathbb{R}^d$ and $T\in[0,1]$,
		\begin{align}\label{acb'}
			\frac{1}{T}\int_{0}^{T}	|b'(t,y)-\bar{b}'(y)|\mathrm{d}t
			+\frac{1}{T}\int_{0}^{T}	|b''(t,y)-\bar{b}''(y)|\mathrm{d}t
			\leq \phi_2(T)(1+|y|).
		\end{align}
	\end{itemize}
	\begin{remark}
		(i) Combining $\textbf{(A3)}$ and $\textbf{(A4)}$, we can verify that the related averaged functions satisfy the assumption $\textbf{(A3)}$, i.e. 
		there exist positive constants $K_7$ and $K_8$ such that for all $y_1, y_2\in\mathbb{R}^d$,
		\begin{align}\label{bb'1}
			|\bar{b}'(y_1)-\bar{b}'(y_2)|+|\bar{b}''(y_1)-\bar{b}''(y_2)|\leq K_7|y_1-y_2|,
		\end{align}
		and
		\begin{align}\label{bb'2}
			|\bar{b}'(y_1)|+|\bar{b}''(y_1)|\leq K_8(1+|y_1|).
		\end{align}	
		(ii) For notational simplicity, denote $K = \max \{K_i\}$ for $i=1,\cdots,8$. 
	\end{remark}

	\subsection{$\Gamma$-convergence}
	We now recall the following fundamental definitions and results of $\Gamma$-convergence. For a comprehensive treatment of the underlying theory, we refer the reader to \cite{braides2002gamma, dal2012introduction}.
	Let $(S, d)$ be a Banach space equipped with a suitable topology $d$, or more generally, a metric space.
	\begin{definition}[$\Gamma$-Convergence \cite{braides2002gamma}]\label{gamma}
		A family of functionals $F_\epsilon: S \to \mathbb{R} \cup \{+\infty\}$ is said to \textit{$\Gamma$-converge} to the limit functional $F_0: S \to \mathbb{R} \cup \{+\infty\}$, as $\epsilon \to 0$, denoted by $F_\epsilon \xrightarrow{\Gamma} F_0$, 
		if the following two conditions are satisfied for every $u \in S$:
		\begin{enumerate}
			\item Liminf inequality: For every sequence $\{u_\epsilon\} \subset S$ such that $u_\epsilon \rightharpoonup u$ weakly in $S$, it holds that
			\begin{equation*}
				\liminf_{\epsilon \to 0} F_\epsilon(u_\epsilon) \ge F_0(u).
			\end{equation*}
			\item Limsup inequality (Recovery sequence): There exists a sequence $\{v_\epsilon\} \subset S$ such that $v_\epsilon \rightharpoonup u$ weakly in $S$ and
			\begin{equation*}
				\limsup_{\epsilon \to 0} F_\epsilon(v_\epsilon) \le F_0(u).
			\end{equation*}
		\end{enumerate}
	\end{definition}
	
	\begin{definition}[Coerciveness Conditions \cite{braides2002gamma}]\label{coerciveness}
		A functional $F: S \to \mathbb{R} \cup \{+\infty\}$ is said to be \textit{coercive} if for all $t \in \mathbb{R}$, the sublevel set $\{u \in S \mid F(u) \le t\}$ is precompact in $S$. 
		
		A functional $F: S \to \mathbb{R} \cup \{+\infty\}$ is \textit{mildly coercive} if there exists a non-empty compact set $K \subset S$ such that
		\begin{equation*}
			\inf_{u \in X} F(u) = \inf_{u \in K} F(u).
		\end{equation*}
		
		Furthermore, a family of functionals $\{F_\epsilon\}_{\epsilon > 0}$ is \textit{equi-mildly coercive} as $\epsilon \to 0$ if there exists a common non-empty compact set $K \subset S$ such that $\inf_{u \in S} F_\epsilon(u) = \inf_{u \in K} F_\epsilon(u)$ for all sufficiently small $\epsilon > 0$.
	\end{definition}
	
	\begin{remark}
		If a functional $F$ is coercive and not identically $+\infty$, it is inherently mildly coercive. We can simply take $K$ as the closure of a non-empty sublevel set $\{F \le t\}$. An intermediate requirement often encountered is that there exists at least one $t \in \mathbb{R}$ such that the set $\{F \le t\}$ is both non-empty and precompact.
	\end{remark}
	
	\begin{lemma}[\cite{braides2002gamma}]\label{min_convergence}
		Let $(S, d)$ be a metric space, and let $\{F_\epsilon\}_{\epsilon > 0}$ be a family of equi-mildly coercive functionals on $S$. If $F_\epsilon \xrightarrow{\Gamma} F_0$ as $\epsilon \to 0$, then the limit functional $F_0$ attains its minimum on $S$, and the minimal values converge
		\begin{equation*}
			\min_{u \in S} F_0(u) = \lim_{\epsilon \to 0} \inf_{u \in S} F_\epsilon(u).
		\end{equation*}
		Moreover, if $\{u_\epsilon\}_{\epsilon > 0} \subset S$ is a precompact sequence such that 
		$$
		\lim_{\epsilon \to 0} F_\epsilon(u_\epsilon) = \lim_{\epsilon \to 0} \inf_{S} F_\epsilon,
		$$
		then every limit of a subsequence of $\{u_\epsilon\}$ as $\epsilon \to 0$ is a minimum point for $F_0$.
	\end{lemma}
	
	\subsection{Generalized Gronwall-Bellman-Bihari inequality}
	Here, we first state the generalization of Gronwall-Bellman-Bihari inequality, which was proved by Pachpatte \cite{pachpatte1975some}.  
	\begin{lemma}[\cite{pachpatte1975some}]\label{geGronwall}
		Let $u$, $f$, $g$ and $h$ be nonnegative continuous functions defined on $\mathbf{R}_{+}$. Let $w(u)$ be a continuous nondecreasing and submultiplicative function defined on $\mathbf{R}_{+}$ and $w(u)>0$ on $(0,\infty)$. If
		$$
		u(t) \leq u_{0}+g(t) \int_{0}^{t} f(s) u(s) \, \mathrm{d} s+\int_{0}^{t} h(s) w(u(s)) \, \mathrm{d} s,
		$$	
		for all $t\in\mathbf{R}_{+}$, where $u_{0}$ is a positive constant, then for $0 \leq t \leq t_{1}$,	
		$$
		u(t) \leq a(t) G^{-1}\left[G\left(u_{0}\right)+\int_{0}^{t} h(s) w(a(s)) \, \mathrm{d} s\right],
		$$
		where
		$$
		a(t)=1+g(t) \int_{0}^{t} f(s) \exp \left(\int_{s}^{t} g(\sigma) f(\sigma) \, \mathrm{d} \sigma\right) \mathrm{d} s,
		$$
		for $t\in\mathbf{R}_{+}$ and
		$$
		G(r)=\int_{r_{0}}^{r} \frac{\mathrm{d} s}{w(s)}, \quad r>0, r_{0}>0,
		$$	
		and $G^{-1}$ is the inverse function of $G$, and $t_{1}\in\mathbf{R}_{+}$ is chosen so that
		$$
		G\left(u_{0}\right)+\int_{0}^{t} h(s) w(a(s)) \, \mathrm{d} s \in \operatorname{Dom}\left(G^{-1}\right),
		$$	
		for all $t\in\mathbf{R}_{+}$ lying in the interval $0 \leq t \leq t_{1}$.
	\end{lemma}
	In particular, we require a specific case of the Lemma presented above.
	\begin{proposition}\label{geGronwall1}
		Let $u$ be a nonnegative continuous function defined on $\mathbf{R}_{+}$. If
		$$
		u(t) \leq u_{0}+C_1\int_{0}^{t} u(s) \, \mathrm{d} s+C_2\int_{0}^{t} u^2(s) \, \mathrm{d} s,
		$$	
		for all $t\in\mathbf{R}_{+}$, where $u_{0}$, $C_1$ and $C_2$ are positive constants, then for $0 \leq t \leq t_{1}=\frac{\ln(1+\frac{2C_1}{u_0C_2})}{2C_1}$,	
		$$
		u(t) \leq \frac{2C_1u_0e^{C_1t}}{2C_1+C_2u_0-C_2u_0e^{2C_1t}}.
		$$
	\end{proposition}
	 \begin{proof}
		 	Applying Lemma \ref{geGronwall} with $g\equiv C_1$, $f\equiv 1$, $h\equiv C_2$ and $w(u)=u^2$, we have 
		 	$$
		 	a(t)=1+C_1 \int_{0}^{t}\exp \left(\int_{s}^{t} C_1 \, \mathrm{d} \sigma\right) \mathrm{d} s
		 	=e^{C_1 t}.
		 	$$
		 	Note that the result is independent of the choice of $r_0$, take $r_0=1$, then it is not hard to obtain 
		 	$$
		 	G(r)=\int_{1}^{r} \frac{1}{s^2}\mathrm{d}s =1-\frac{1}{r},\quad \text{with $r>0$},
		 	$$
		 	and its inverse
		 	$$
		 	G^{-1}(y)=\frac{1}{1-y}, \quad \text{with $y<1$}.
		 	$$
		 	Hence, there exists $t_1=\frac{\ln(1+\frac{2C_1}{u_0C_2})}{2C_1}$ such that for all $0\leq t\leq t_1$, it holds that
		 	\begin{align*}
			 		u(t) &\leq a(t) G^{-1}\left[G\left(u_{0}\right)+C_2\int_{0}^{t} a^2(s) \, \mathrm{d} s\right]
			 		\leq e^{C_1 t} G^{-1}\left[1-\frac{1}{u_{0}}+C_2\int_{0}^{t} e^{2C_1 s} \, \mathrm{d} s\right]\\
			 		&=e^{C_1 t} G^{-1}\left[1-\frac{1}{u_{0}}+\frac{C_2[e^{2C_1t} -1]}{2C_1}\right]
			 		= e^{C_1 t}\frac{1}{\frac{1}{u_{0}}-\frac{C_2[e^{2C_1t} -1]}{2C_1}}
			 		=\frac{2C_1u_0e^{C_1t}}{2C_1+C_2u_0-C_2u_0e^{2C_1t}}.
			 	\end{align*}	
		 	Indeed, $t_1$ can be solved by
		 	$$
		 	G\left(u_{0}\right)+C_2\int_{0}^{t} a^2(s) \, \mathrm{d} s 
		 	=1-\frac{1}{u_{0}}+ C_2\int_{0}^{t_1} e^{2C_1 s} \, \mathrm{d} s
		 	=1-\frac{1}{u_{0}}+\frac{C_2[e^{2C_1t} -1]}{2C_1}<1.
		 	$$
		 	The proof is complete.
		 \end{proof}
	
	\section{Proximity and convergence of the most probable transition paths}\label{PC}
In this section, we study the relationship between the most probable transition paths of the multiscale system and those of the averaged system under the supremum norm, based on Euler-Lagrange equations. 
We first derive several estimates for these paths. We then establish a proximity result and finally prove a convergence result.

	%

	\subsection{Some estimates of most probable transition paths}
	In this subsection, we will prove some estimates of MPTPs. Specially, we prove that  the processes $Z^\epsilon$, $\bar{Z}^\epsilon$, $\varPhi^\epsilon$ and $\bar{\varPhi}$ are bounded with fixed $0<\epsilon<1$ on the time horizon $[0,1]$. 
	From now on, to keep notation simple, we focus only on the one-dimensional case, but the higher-dimensional case follows similarly.
	\begin{lemma}\label{Zbounded}
		Assume that $\textbf{(A1)}$ and $\textbf{(A3)}$ hold.
		For any fixed $\epsilon$ and all $\tau\in[0,1]$, 
		the MPTP $Z^\epsilon$ is bounded by a positive constant $C$ depending on $|x_0|$, $|\dot{Z}_0|$, $K$ and $\sigma$, that is 
		$$\sup_{0\leq \tau\leq 1}|Z^\epsilon_\tau|\leq \frac{2C_1u_0e^{C_1}}{2C_1+C_2u_0-C_2u_0e^{2C_1}}:=C,$$ 
		where 
		$u_0=|x_0| + |\dot{Z}_0| +\frac{ K^2}{2} +\frac{K}{2}+K$, 
		$C_1=\frac{ K}{2}+K$, 
		$C_2=2K^2$.
	\end{lemma}
	\begin{proof}
		Without loss of generality, we only prove the case of $\epsilon=1$. Then Equ.\eqref{EL1} becomes
		\begin{align*}
			\ddot{Z}_t =  b(t,Z_t)b'(t,Z_t)
			+\frac{1}{2} b''(t,Z_t)
			+\partial_t b(t,Z_t).
		\end{align*}
		For $s\in[0,1]$, integrating both sides to the above equation with respect to $t$, we obtain
		\begin{align*}
			\dot{Z}_s
			= \dot{Z}_0 
			+ \int_0^s b(t,Z_t)b'(t,Z_t)
			\mathrm{d}t
			+\frac{1}{2} \int_0^s b''(t,Z_t)\mathrm{d}t
			+\int_0^s \partial_t b(t,Z_t)\mathrm{d}t.
		\end{align*}
		For $\tau\in[0,1]$, integrating again yields 
		\begin{align*}
			Z_\tau
			= Z_0 +\tau\dot{Z}_0 
			+ \int_0^\tau\int_0^s b(t,Z_t)b'(t,Z_t)
			\mathrm{d}t\mathrm{d}s
			+\frac{1}{2} \int_0^\tau\int_0^s b''(t,Z_t)
			\mathrm{d}t\mathrm{d}s
			+\int_0^\tau\int_0^s \partial_t b(t,Z_t)
			\mathrm{d}t\mathrm{d}s.
		\end{align*}
		Therefore, by growth conditions \eqref{b2}, \eqref{b'2} and $(a+b)^2 \leq 2(a^2+b^2)$, we have
		\begin{align*}
			|Z_\tau|
			&\leq |x_0| +\tau|\dot{Z}_0| 
			+ \left|\int_0^\tau\int_0^s b(t,Z_t)b'(t,Z_t)
			\mathrm{d}t\mathrm{d}s\right|
			+\frac{1}{2} \left|\int_0^\tau\int_0^s b''(t,Z_t)
			\mathrm{d}t\mathrm{d}s\right|
			+\left|\int_0^\tau\int_0^s \partial_t b(t,Z_t)
			\mathrm{d}t\mathrm{d}s\right|\\
			&\leq |x_0| +\tau|\dot{Z}_0| 
			+ \int_0^\tau\int_0^s |b(t,Z_t)||b'(t,Z_t)|
			\mathrm{d}t\mathrm{d}s
			+\frac{1}{2} \int_0^\tau\int_0^s |b''(t,Z_t)|
			\mathrm{d}t\mathrm{d}s
			+\int_0^\tau\int_0^s |\partial_t b(t,Z_t)|
			\mathrm{d}t\mathrm{d}s\\
			&\leq |x_0| +\tau|\dot{Z}_0| 
			+\int_0^\tau\int_0^s  K_2(t)(1+|Z_t|) K_6(t)(1+|Z_t|)
			\mathrm{d}t\mathrm{d}s
			+\frac{3}{2} \int_0^\tau\int_0^s K_6(t)(1+|Z_t|)
			\mathrm{d}t\mathrm{d}s\\
			&\leq |x_0| +\tau|\dot{Z}_0| 
			+ \int_0^\tau (\tau-t) K_2(t)K_6(t)(1+|Z_t|)^2 \mathrm{d}t
			+\frac{3}{2} \int_0^\tau (\tau-t) K_6(t)(1+|Z_t|)\mathrm{d}t\\
			&\leq |x_0| +\tau|\dot{Z}_0| 
			+ \int_0^\tau  K_2(t)K_6(t)(1+|Z_t|)^2 \mathrm{d}t
			+\frac{3}{2} \int_0^\tau K_6(t)(1+|Z_t|)\mathrm{d}t\\
			&\leq |x_0| +\tau|\dot{Z}_0| 
			+ 2\int_0^\tau  K_2(t)K_6(t)(1+|Z_t|^2) \mathrm{d}t
			+\frac{3}{2}\int_0^\tau K_6(t)(1+|Z_t|)\mathrm{d}t\\
			&\leq |x_0| + |\dot{Z}_0| 
			+2K^2 \int_0^\tau |Z_t|^2 \mathrm{d}t
			+\frac{3K}{2}\int_0^\tau  |Z_t| \mathrm{d}t 
			+2K^2 +\frac{3K}{2}.
		\end{align*}
		Applying Proposition \ref{geGronwall1} to $u=|Z|$, 
		$u_0=|x_0| + |\dot{Z}_0| +2K^2 +\frac{3K}{2}$, 
		$C_1=\frac{3K}{2}$, 
		$C_2=2K^2$ and choosing a suitable $K$ satisfying $t_1=\frac{\ln(1+\frac{2C_1}{u_0C_2})}{2C_1}\geq 1$,  we get for all $\tau\in[0,1]$,
		\begin{align*}
			\sup_{0\leq \tau\leq 1}	|Z_\tau|\leq \frac{2C_1u_0e^{C_1\tau}}{2C_1+C_2u_0-C_2u_0e^{2C_1\tau}} 
			\leq \frac{2C_1u_0e^{C_1}}{2C_1+C_2u_0-C_2u_0e^{2C_1}}:=C,
		\end{align*}
		which complete the proof.
	\end{proof}
	\begin{lemma}\label{Psibounded}
		Assume that $\textbf{(A1)}$ and $\textbf{(A3)}$ hold.
		For any fixed $\epsilon$ and all $\tau\in[0,1]$, the MPTP $\varPhi^\epsilon$ is bounded by a positive constant $C$ depending on $|x_0|$, $|\dot{\varPhi}_0|$, $K$ and $\sigma$, that is 
		$$\sup_{0\leq \tau\leq 1}|\varPhi^\epsilon_\tau|\leq \frac{2C_1u_0e^{C_1}}{2C_1+C_2u_0-C_2u_0e^{2C_1}}:=C,$$ 
		where $u_0=|x_0| + |\dot{\varPhi}_0| +\frac{ K^2}{2} +\frac{3K}{2}$, 
		$C_1=\frac{3K}{2}$ and
		$C_2=2K^2$.
	\end{lemma}
	\begin{proof}
		The path $\varPhi^\epsilon$ is the same as $Z^\epsilon$ when $\epsilon=1$, the proof is omitted here.
	\end{proof}
	
	\begin{lemma}\label{Zbbounded}
		Assume that $\textbf{(A1)-(A4)}$ hold.
		For any fixed $\epsilon$ and all $\tau\in[0,1]$, the MPTP $\bar{Z}^\epsilon$ is bounded by a positive constant $\bar{C}$ depending on $|x_0|$, $|\dot{\bar{Z}}_0|$, $K$ and $\bar{\sigma}$, that is
		$$\sup_{0\leq \tau\leq 1}|\bar{Z}^\epsilon_\tau|\leq \frac{2\bar{C}_1\bar{u}_0e^{\bar{C}_1}}{2\bar{C}_1+\bar{C}_2\bar{u}_0-\bar{C}_2\bar{u}_0e^{2\bar{C}_1}}:=\bar{C},$$ 
		where $\bar{u}_0=|x_0| + |\dot{\bar{Z}}_0| +2K^2 +\frac{  K}{2}$, $\bar{C}_1=\frac{  K}{2}$ and $\bar{C}_2=2K^2$.
	\end{lemma}
	\begin{proof}
		By Equ.\eqref{EL2}, the triangle inequality and Lemma \ref{Zbounded}, it can be deduced that
		\begin{align*}
			|\bar{Z}_\tau|
			&\leq |x_0| +\tau|\dot{\bar{Z}}_0| 
			+ \left|\int_0^\tau\int_0^s \bar{b}(\bar{Z}_t)\bar{b}'(\bar{Z}_t)
			\mathrm{d}t\mathrm{d}s\right|
			+\frac{1}{2} \left|\int_0^\tau\int_0^s \bar{b}''(\bar{Z}_t)
			\mathrm{d}t\mathrm{d}s\right|\\
			&\leq |x_0| +\tau|\dot{\bar{Z}}_0| 
			+ 2\int_0^\tau  K_2(t)K_6(t)(1+|\bar{Z}_t|^2) \mathrm{d}t
			+\frac{1}{2} \int_0^\tau K_6(t)(1+|\bar{Z}_t|)\mathrm{d}t\qquad\quad\\
			&\leq |x_0| + |\dot{\bar{Z}}_0| 
			+2K^2 \int_0^\tau |\bar{Z}_t|^2 \mathrm{d}t
			+\frac{  K}{2}\int_0^\tau  |\bar{Z}_t| \mathrm{d}t 
			+2K^2 +\frac{  K}{2}.
		\end{align*}
		Since all $K_i$ have the same upper bound $K>0$ for $i=1,\cdots,8$, the result can be easliy obtained with growth conditions \eqref{bb2}, \eqref{bb'2} and the similiar procedure in Lemma \ref{Zbounded}. We omit the details.
	\end{proof}
	\begin{lemma}\label{Psibbounded}
		Assume that $\textbf{(A1)-(A4)}$ hold.
		For all $\tau\in[0,1]$, the MPTP $\bar{\varPhi}$ is bounded by a positive constant $\bar{C}$ depending on $|x_0|$, $|\dot{\bar{\varPhi}}_0|$, $K$ and $\bar{\sigma}$, that is,
		$$\sup_{0\leq \tau\leq 1}|\bar{\varPhi}_\tau|\leq \frac{2\bar{C}_1\bar{u}_0e^{\bar{C}_1}}{2\bar{C}_1+\bar{C}_2\bar{u}_0-\bar{C}_2\bar{u}_0e^{2\bar{C}_1}}:=\bar{C},$$ 
		where $\bar{u}_0=|x_0| + |\dot{\bar{\varPhi}}_0| +2K^2 +\frac{  K}{2}$, $\bar{C}_1=\frac{  K}{2}$ and $\bar{C}_2=2K^2$.
	\end{lemma}
	\begin{proof}
		The path $\bar{\varPhi}$ is the same as $\bar{Z}^\epsilon$ when $\epsilon=1$, the proof is omitted here.
	\end{proof}
	
	\begin{remark}
		From Lemma \ref{Zbounded}, Lemma \ref{Psibounded}, Lemma \ref{Zbbounded} and Lemma \ref{Psibbounded}, the existence and uniqueness of solutions to Equ.\eqref{EL1}, Equ.\eqref{EL11}, Equ.\eqref{EL2} and Equ.\eqref{EL3} can be guaranteed.
	\end{remark}
	\subsection{Proximity between the most probable transition paths}
	We now consider the connection between the processes $Z^\epsilon_t$ and $\bar{Z}^\epsilon_t$. Since both processes depend on the parameter $\epsilon$, we shall be interested in their pointwise proximity, which is the first result as follows. 
	
	\begin{theorem}[Proximity Result]\label{main}
		Assume that $\textbf{(A1)-(A4)}$ hold. If $Z^\epsilon$ and $\bar{Z}^\epsilon$ are the most probable transition path for SDE \eqref{SDE1} and averaged SDE \eqref{SDE2}, which are characterized by Euler-Langrange equations \eqref{EL1} and \eqref{EL2}, then we can derive the proximity between the paths in $C([0,1];\mathbb{R}^d)$, more specifically, there exist constants $C_p$ and $C_q$ such that 
		$$
		\sup_{\tau \in (0, 1)}|Z^\epsilon_\tau - \bar{Z}^\epsilon_\tau|
		\leq C_p\epsilon e^{C_q\epsilon^2},
		$$
		where $C_p$ and $C_q$ depend on $C$, $\bar{C}$, $K$ and the upper bound of $\phi_1$, $\phi_2$.
		%
	\end{theorem}
	\begin{proof}
		We only prove the one-dimensional case. According to the Euler-Lagrange equations \eqref{EL11} and \eqref{EL2}, we have
		\begin{align*}
			\ddot{Z}^\epsilon_t - \ddot{\bar{Z}}^\epsilon_t 
			=\epsilon^2 \left[b(t,Z^{\epsilon}_t)b'(t,Z^{\epsilon}_t)
			- \bar{b}(\bar{Z}^{\epsilon}_t)\bar{b}'(\bar{Z}^{\epsilon}_t)\right]
			+\frac{  \epsilon^2}{2} b''(t,Z^{\epsilon}_t)
			-\frac{  \epsilon^2}{2} \bar{b}''(\bar{Z}^{\epsilon}_t)
			+\epsilon\partial_t b(t,Z^{\epsilon}_t).
		\end{align*}
		For $s\in[0,1]$, integrating both sides with respect to $t$, we obtain
		\begin{align*}
			\dot{Z}^\epsilon_s - \dot{\bar{Z}}^\epsilon_s
			=&(\dot{Z}^\epsilon_0 - \dot{\bar{Z}}^\epsilon_0 ) 
			+\epsilon^2 \int_0^s\left[b(t,Z^{\epsilon}_t)b'(t,Z^{\epsilon}_t)
			- \bar{b}(\bar{Z}^{\epsilon}_t)\bar{b}'(\bar{Z}^{\epsilon}_t)\right]
			\mathrm{d}t\nonumber\qquad\qquad\qquad\qquad\quad\\
			&+\frac{\epsilon^2}{2} \int_0^s \left[  b''(t,Z^{\epsilon}_t)
			-  \bar{b}''(\bar{Z}^{\epsilon}_t)\right]\mathrm{d}t
			+\epsilon  \int_0^s \partial_t b(t,Z^{\epsilon}_t)
			\mathrm{d}t.
		\end{align*}
		For $\tau\in[0,1]$, integrating again yields 
		\begin{align*}
			Z^\epsilon_\tau - \bar{Z}^\epsilon_\tau
			=&\left(Z^\epsilon_0 - \bar{Z}^\epsilon_0\right)
			+\tau(\dot{Z}^\epsilon_0 - \dot{\bar{Z}}^\epsilon_0)
			+\epsilon^2 \int_0^\tau\int_0^s\left[b(t,Z^{\epsilon}_t)b'(t,Z^{\epsilon}_t)
			- \bar{b}(\bar{Z}^{\epsilon}_t)\bar{b}'(\bar{Z}^{\epsilon}_t)\right]
			\mathrm{d}t\mathrm{d}s\\
			&+\frac{\epsilon^2}{2} \int_0^\tau\int_0^s \left[  b''(t,Z^{\epsilon}_t)
			-  \bar{b}''(\bar{Z}^{\epsilon}_t)\right]\mathrm{d}t\mathrm{d}s
			+\epsilon  \int_0^\tau\int_0^s \partial_t b(t,Z^{\epsilon}_t)\mathrm{d}t\mathrm{d}s\\
			=&\left(Z^\epsilon_0 - \bar{Z}^\epsilon_0\right)
			+\tau(\dot{Z}^\epsilon_0 - \dot{\bar{Z}}^\epsilon_0)+I_{1}+I_{2}+I_{3}+I_{4}+I_{5}+I_{6}+I_{7},
		\end{align*}
		where 
		\begin{align*}
			&I_{1}=\epsilon^2 \int_0^\tau\int_0^s 
			b'(t,Z^{\epsilon}_t)\left[b(t,Z^{\epsilon}_t)- \bar{b}(Z^{\epsilon}_t)\right]
			\mathrm{d}t\mathrm{d}s,\quad
			I_{2}=\epsilon^2 \int_0^\tau\int_0^s 
			\bar{b}(Z^{\epsilon}_t)\left[b'(t,Z^{\epsilon}_t)-\bar{b}'(Z^{\epsilon}_t)\right]
			\mathrm{d}t\mathrm{d}s,\\
			&I_{3}=\epsilon^2 \int_0^\tau\int_0^s 
			\bar{b}'(Z^{\epsilon}_t)\left[\bar{b}(Z^{\epsilon}_t)-\bar{b}(\bar{Z}^{\epsilon}_t)\right]
			\mathrm{d}t\mathrm{d}s,\quad\quad\,\,\,
			I_{4}=\epsilon^2 \int_0^\tau\int_0^s 
			\bar{b}(\bar{Z}^{\epsilon}_t)\left[\bar{b}'(Z^{\epsilon}_t)-\bar{b}'(\bar{Z}^{\epsilon}_t)\right]
			\mathrm{d}t\mathrm{d}s,\\
			&I_{5}=\frac{\epsilon^2  }{2} \int_0^\tau\int_0^s 
			b''(t,Z^{\epsilon}_t)\mathrm{d}t\mathrm{d}s,\quad
			I_{6}=-\frac{\epsilon^2  }{2} \int_0^\tau\int_0^s 
			\bar{b}''(\bar{Z}^{\epsilon}_t)\mathrm{d}t\mathrm{d}s,\quad
			I_{7}=\epsilon  \int_0^\tau\int_0^s \partial_t b(t,Z^{\epsilon}_t)\mathrm{d}t\mathrm{d}s.
		\end{align*}
		For the term $|I_{1}|$, by \eqref{ac}, \eqref{b'2} and Lemma \ref{Zbounded}, we obtain
		\begin{align*}
			|I_{1}|&
			\leq \epsilon^2 \int_0^\tau\!\int_0^s \!
			K(1+|Z^{\epsilon}_t|)|b(t,Z^{\epsilon}_t)- \bar{b}(Z^{\epsilon}_t)| 
			\mathrm{d}t\mathrm{d}s
			\leq \epsilon^2(1+C)K \int_0^\tau\! \frac{1}{s}\!\int_0^s 
			|b(t,Z^{\epsilon}_t)- \bar{b}(Z^{\epsilon}_t)| 
			\mathrm{d}t\mathrm{d}s\\
			&
			\leq \epsilon^2(1+C)K \int_0^\tau \phi_1(s)(1+\sup_{0\leq s\leq 1}|Z^{\epsilon}_s|)\mathrm{d}s
			\leq \epsilon^2(1+C)^2K \int_0^\tau \phi_1(s) \mathrm{d}s:=\epsilon^2 C_{11},
		\end{align*}
		where $C_{11}$ is a positive constant depedending on $C$, $K$ and the upper bound of $\phi_1$.
		
		For the term $|I_{2}|$, by \eqref{bb2}, \eqref{acb'} and Lemma \ref{Zbounded}, we obtain
		\begin{align*}
			|I_{2}|&  
			\leq \epsilon^2 \!\int_0^\tau\!\int_0^s \!
			K(1+|Z^{\epsilon}_t|)|b'(t,Z^{\epsilon}_t)- \bar{b}'(Z^{\epsilon}_t)| 
			\mathrm{d}t\mathrm{d}s
			\leq \epsilon^2(1+C)K\! \int_0^\tau\! \frac{1}{s}\int_0^s\! 
			|b'(t,Z^{\epsilon}_t)- \bar{b}'(Z^{\epsilon}_t)| 
			\mathrm{d}t\mathrm{d}s\\
			&
			\leq \epsilon^2(1+C)K \int_0^\tau \phi_2(s)(1+\sup_{0\leq s\leq 1}|Z^{\epsilon}_s|)\mathrm{d}s
			\leq \epsilon^2(1+C)^2K \int_0^\tau \phi_2(s) \mathrm{d}s:=\epsilon^2 C_{21},
		\end{align*}
		where $C_{21}$ is a positive constant depedending on $C$, $K$ and the upper bound of $\phi_2$.
		
		For the term $|I_{3}|$, by \eqref{bb1}, \eqref{bb'2} and Lemma \ref{Zbounded}, we obtain
		\begin{align*}
			|I_{3}|&
			\leq \epsilon^2 \int_0^\tau\int_0^s 
			K(1+|Z^{\epsilon}_t|)K|Z^{\epsilon}_t- \bar{Z}^{\epsilon}_t| 
			\mathrm{d}t\mathrm{d}s
			\leq \epsilon^2(1+C)K^2 \int_0^\tau\int_0^s 
			|Z^{\epsilon}_t- \bar{Z}^{\epsilon}_t| 
			\mathrm{d}t\mathrm{d}s\qquad\quad\qquad\\
			&
			=\epsilon^2(1+C)K^2 \int_0^\tau (\tau-t) 
			|Z^{\epsilon}_t- \bar{Z}^{\epsilon}_t| 
			\mathrm{d}t
			\leq\epsilon^2(1+C)K^2 \int_0^\tau 
			|Z^{\epsilon}_t- \bar{Z}^{\epsilon}_t| 
			\mathrm{d}t.
		\end{align*}
		
		For the term $|I_{4}|$, by \eqref{bb2}, \eqref{bb'1} and Lemma \ref{Zbbounded}, we obtain
		\begin{align*}
			|I_{4}|&
			\leq \epsilon^2 \int_0^\tau\int_0^s 
			K(1+|\bar{Z}^{\epsilon}_t|)K|Z^{\epsilon}_t- \bar{Z}^{\epsilon}_t| 
			\mathrm{d}t\mathrm{d}s
			\leq \epsilon^2(1+\bar{C})K^2 \int_0^\tau\int_0^s 
			|Z^{\epsilon}_t- \bar{Z}^{\epsilon}_t| 
			\mathrm{d}t\mathrm{d}s\qquad\qquad\\
			&
			=\epsilon^2(1+\bar{C})K^2 \int_0^\tau (\tau-t) 
			|Z^{\epsilon}_t- \bar{Z}^{\epsilon}_t| 
			\mathrm{d}t
			\leq\epsilon^2(1+\bar{C})K^2 \int_0^\tau 
			|Z^{\epsilon}_t- \bar{Z}^{\epsilon}_t| 
			\mathrm{d}t.
		\end{align*}
		
		For the term $|I_{5}|$, $|I_{6}|$ and $|I_{7}|$, by \eqref{b'2}, \eqref{bb'2}, Lemma \ref{Zbounded} and Lemma \ref{Zbbounded}, we obtain
		\begin{align*}
			&|I_{5}|+|I_{6}|+|I_{7}|
			\leq \frac{\epsilon^2  }{2} \int_0^\tau 
			|b''(t,Z^{\epsilon}_t)|\mathrm{d}t 
			+\frac{\epsilon^2  }{2} \int_0^\tau 
			|\bar{b}''(\bar{Z}^{\epsilon}_t)|\mathrm{d}t 
			+\epsilon  \int_0^\tau |\partial_t b(t,Z^{\epsilon}_t)|\mathrm{d}t \qquad\quad
			\\
			\leq& \frac{\epsilon^2  K}{2} \int_0^\tau 
			(1+|Z^{\epsilon}_t|)\mathrm{d}t 
			+\frac{\epsilon^2  K}{2} \int_0^\tau 
			(1+|\bar{Z}^{\epsilon}_t|)\mathrm{d}t 
			+\epsilon K  \int_0^\tau (1+|Z^{\epsilon}_t|)\mathrm{d}t 
			\\
			=&\left(\frac{\epsilon^2  K}{2}+\epsilon K\right)
			\int_0^\tau |Z^{\epsilon}_t|\mathrm{d}t 
			+\frac{\epsilon^2  K}{2}\int_0^\tau 
			|\bar{Z}^{\epsilon}_t|\mathrm{d}t 
			+\frac{\epsilon^2  K}{2}+\frac{\epsilon^2  K}{2}+\epsilon K\\
			\leq& \left(\frac{\epsilon^2  K}{2}+\epsilon K\right)C
			+\frac{\epsilon^2  K}{2}\bar{C}
			+\epsilon^2  K+\epsilon K.
		\end{align*}
		
		Hence, together with the triangle inequality, we have
		\begin{align*}
			\quad&|Z^\epsilon_\tau - \bar{Z}^\epsilon_\tau|\leq 
			|Z^\epsilon_0 - \bar{Z}^\epsilon_0|
			+|\dot{Z}^\epsilon_0 - \dot{\bar{Z}}^\epsilon_0|
			+|I_{1}|+|I_{2}|+|I_{3}|+|I_{4}|+|I_{5}|+|I_{6}|+|I_{7}|\\
			\leq& 
			|Z^\epsilon_0 - \bar{Z}^\epsilon_0|
			+|\dot{Z}^\epsilon_0 - \dot{\bar{Z}}^\epsilon_0|
			+\epsilon^2 C_{11}+\epsilon^2 C_{12}
			+\left[\epsilon^2(1+C)K^2 
			+\epsilon^2(1+\bar{C})K^2 
			\right]
			\int_0^\tau 
			|Z^{\epsilon}_t- \bar{Z}^{\epsilon}_t| 
			\mathrm{d}t 
			\\
			& +
			\left(\frac{\epsilon^2  K}{2}+\epsilon K\right)C
			+\frac{\epsilon^2  K}{2}\bar{C}
			+\epsilon^2  K+\epsilon K,
		\end{align*}
		with the aid of Gronwall inequality, it follows that
		\begin{align*}
			\sup_{\tau \in (0, 1)}|Z^\epsilon_\tau - \bar{Z}^\epsilon_\tau|&\leq 
			\left[ |\dot{Z}^\epsilon_0 - \dot{\bar{Z}}^\epsilon_0|
			+\epsilon^2 C_{11}+\epsilon^2 C_{12}
			+\left(\frac{\epsilon^2  K}{2}+\epsilon K\right)C
			+\frac{\epsilon^2  K}{2}\bar{C}
			+\epsilon^2  K+\epsilon K\right]\\
			&\quad \times \exp\{\epsilon^2(1+C)K^2 
			+\epsilon^2(1+\bar{C})K^2 \}
			:=C_p\epsilon e^{C_q\epsilon^2},
		\end{align*}
		where the constants $C_p$ and $C_q$ depend on $C$, $\bar{C}$, $K$ and the upper bound of $\phi_1$, $\phi_2$. 
	\end{proof}
	
	\subsection{Convergence of the most probable transition paths}
	Next we will prove the convergence of the MPTPs $\varPhi^\epsilon_{t}$
	under the supremum norm, which is the second result. 
	\begin{theorem}[Convergence Result]\label{mainn}
		Assume that $\textbf{(A1)-(A4)}$ hold. If $\varPhi^\epsilon$ and $\bar{\varPhi}$ are the most probable transition paths for SDE \eqref{SDE11} and averaged SDE \eqref{SDE3}, which are characterized by Euler-Langrange equations \eqref{EL11} and \eqref{EL3}, furthermore,
		if  $b(0,x_0)=0$ and the time average of the drift $b(t,y)$ vanishes in the long-time limit for every $y\in \mathbb{R}^d$, 
		then we can derive the convergence between the paths in $C([0,1];\mathbb{R}^d)$, more specifically, we have
		$$\lim\limits_{\epsilon \rightarrow 0}\sup_{\tau \in [0, 1]}
		|\varPhi^\epsilon_\tau - \bar{\varPhi}_\tau|=0.$$
	\end{theorem}
	\begin{proof}
		According to the Euler--Lagrange equations \eqref{EL1} and \eqref{EL3}, we have
		\begin{align*}\label{ddotPsi}
			\ddot{\varPhi}^\epsilon_t - \ddot{\bar{\varPhi}}_t 
			=\left[b(t/\epsilon,\varPhi^{\epsilon}_t)b'(t/\epsilon,\varPhi^{\epsilon}_t)
			- \bar{b}(\bar{\varPhi}_t)\bar{b}'(\bar{\varPhi}_t)\right]
			+\frac{1}{2} b''(t/\epsilon,\varPhi^{\epsilon}_t)
			-\frac{1}{2}\bar{b}''(\bar{\varPhi}_t)
			+\frac{1}{\epsilon}\partial_t b(t/\epsilon,\varPhi^{\epsilon}_t).
		\end{align*}
		For $s\in[0,1]$, integrating both sides with respect to $t$, we obtain
		\begin{align*}
			\dot{\varPhi}^\epsilon_s - \dot{\bar{\varPhi}}_s
			=&(\dot{Z}^\epsilon_0 - \dot{\bar{Z}}^\epsilon_0 ) 
			+\int_0^s\left[b(t/\epsilon,\varPhi^{\epsilon}_t)b'(t/\epsilon,\varPhi^{\epsilon}_t)
			- \bar{b}(\bar{\varPhi}_t)\bar{b}'(\bar{\varPhi}_t)\right]
			\mathrm{d}t\nonumber\qquad\qquad\qquad\qquad\qquad\\
			&+\frac{1}{2} \int_0^s \left[  b''(t/\epsilon,\varPhi^{\epsilon}_t)
			-  \bar{b}''(\bar{\varPhi}_t)\right]\mathrm{d}t
			+\frac{1}{\epsilon} \int_0^s \partial_t b(t/\epsilon,\varPhi^{\epsilon}_t)
			\mathrm{d}t.
		\end{align*}
		For $\tau\in[0,1]$, integrating again yields 
		\begin{align*}
			\varPhi^\epsilon_\tau - \bar{\varPhi}_\tau
			=&\left(\varPhi^\epsilon_0 - \bar{\varPhi}_0\right)
			+\tau(\dot{\varPhi}^\epsilon_0 - \dot{\bar{\varPhi}}_0)
			+\int_0^\tau\int_0^s\left[b(t/\epsilon,\varPhi^{\epsilon}_t)b'(t/\epsilon,\varPhi^{\epsilon}_t)
			- \bar{b}(\bar{\varPhi}_t)\bar{b}'(\bar{\varPhi}_t)\right]
			\mathrm{d}t\mathrm{d}s\quad\\
			&+\frac{1}{2} \int_0^\tau\int_0^s \left[  b''(t/\epsilon,\varPhi^{\epsilon}_t)
			-  \bar{b}''(\bar{\varPhi}_t)\right]\mathrm{d}t\mathrm{d}s
			+\frac{1}{\epsilon}\int_0^\tau\int_0^s \partial_t b(t/\epsilon,\varPhi^{\epsilon}_t)\mathrm{d}t\mathrm{d}s\nonumber\\
			=&\left(\varPhi^\epsilon_0 - \bar{\varPhi}_0\right)
			+\tau(\dot{\varPhi}^\epsilon_0 - \dot{\bar{\varPhi}}_0)+J_{1}+J_{2}+J_{3}+J_{4}+J_{5}+J_{6}+J_{7},
		\end{align*}
		where 
		\begin{align*}
			&J_{1}=\int_0^\tau\!\int_0^s \!
			b'(t/\epsilon,\varPhi^{\epsilon}_t)\!\left[b(t/\epsilon,\varPhi^{\epsilon}_t)- \bar{b}(\varPhi^{\epsilon}_t)\right]\!
			\mathrm{d}t\mathrm{d}s,\quad
			J_{2}=\int_0^\tau\!\int_0^s \!
			\bar{b}(\varPhi^{\epsilon}_t)\!\left[b'(t/\epsilon,\varPhi^{\epsilon}_t)-\bar{b}'(\varPhi^{\epsilon}_t)\right]\!
			\mathrm{d}t\mathrm{d}s,\\
			&J_{3}=\int_0^\tau\!\int_0^s \!
			\bar{b}'(\varPhi^{\epsilon}_t)\!\left[\bar{b}(\varPhi^{\epsilon}_t)-\bar{b}(\bar{\varPhi}_t)\right]\!
			\mathrm{d}t\mathrm{d}s,\quad\quad\quad\quad\;\;
			J_{4}=\int_0^\tau\!\int_0^s\! 
			\bar{b}(\bar{\varPhi}_t)\!\left[\bar{b}'(\varPhi^{\epsilon}_t)-\bar{b}'(\bar{\varPhi}_t)\right]\!
			\mathrm{d}t\mathrm{d}s,\\
			&J_{5}=\frac{1}{2} \int_0^\tau\int_0^s 
			\left[b''(t/\epsilon,\varPhi^{\epsilon}_t)-b''(t/\epsilon,\bar{\varPhi}_t)\right]\mathrm{d}t\mathrm{d}s,\quad\;
			J_{6}=\frac{1}{2} \int_0^\tau\int_0^s 
			\left[b''(t/\epsilon,\bar{\varPhi}_t)-\bar{b}''(\bar{\varPhi}_t)\right]\mathrm{d}t\mathrm{d}s,\\
			&J_{7}=\frac{1}{\epsilon}  \int_0^\tau\int_0^s \partial_t b(t/\epsilon,\varPhi^{\epsilon}_t)\mathrm{d}t\mathrm{d}s.
		\end{align*}
		
		For the term $J_{1}$, by exchanging the order of integration, change of variables, \eqref{ac}, \eqref{b'2} and Lemma \ref{Psibounded}, we obtain
		\begin{align*}
			|J_{1}|&
			\leq \int_0^\tau\! \int_0^s \! 
			K(1+|\varPhi^{\epsilon}_t|)|b(t/\epsilon,\varPhi^{\epsilon}_t)- \bar{b}(\varPhi^{\epsilon}_t)| 
			\mathrm{d}t\mathrm{d}s
			\leq (1+C)K\!  \int_0^\tau\!  (\tau-t)
			|b(t/\epsilon,\varPhi^{\epsilon}_t)- \bar{b}(\varPhi^{\epsilon}_t)| 
			\mathrm{d}t\\
			&\leq (1+C)K\! \int_0^\tau\! 
			|b(t/\epsilon,\varPhi^{\epsilon}_t)- \bar{b}(\varPhi^{\epsilon}_t)| 
			\mathrm{d}t
			\leq (1+C)K \frac{\epsilon}{\tau}\int_0^{\frac{\tau}{\epsilon}}
			|b(u,\varPhi^{\epsilon}_{u\epsilon})- \bar{b}(\varPhi^{\epsilon}_{u\epsilon})| 
			\mathrm{d}u\\
			&\leq (1+C)^2 K \phi_1\left(\frac{\tau}{\epsilon}\right).
		\end{align*}
		
		For the term $J_{2}$, by \eqref{bb2}, \eqref{acb'} and Lemma \ref{Psibounded}, we obtain
		\begin{align*}
			|J_{2}|&\leq \int_0^\tau\int_0^s 
			|\bar{b}(\varPhi^{\epsilon}_t)||b'(t/\epsilon,\varPhi^{\epsilon}_t)-\bar{b}'(\varPhi^{\epsilon}_t)|\mathrm{d}t\mathrm{d}s
			\leq \int_0^\tau\int_0^s 
			K(1+|\varPhi^{\epsilon}_t|)|b'(t/\epsilon,\varPhi^{\epsilon}_t)- \bar{b}'(\varPhi^{\epsilon}_t)| 
			\mathrm{d}t\mathrm{d}s\\
			&\leq (1+C)K \frac{\epsilon}{\tau}\int_0^{\frac{\tau}{\epsilon}}
			|b'(u,\varPhi^{\epsilon}_{u\epsilon})- \bar{b}'(\Phi^{\epsilon}_{u\epsilon})| 
			\mathrm{d}u
			\leq (1+C)^2 K \phi_2\left(\frac{\tau}{\epsilon}\right).
		\end{align*}
		
		For the term $J_{3}$, by \eqref{bb1}, \eqref{bb'2} and Lemma \ref{Psibounded}, we obtain
		\begin{align*}
			|J_{3}|
			&\leq\int_0^\tau\int_0^s 
			K(1+|\varPhi^{\epsilon}_t|)K|\varPhi^{\epsilon}_t- \bar{\varPhi}_t| 
			\mathrm{d}t\mathrm{d}s
			\leq (1+C)K^2 \int_0^\tau\int_0^s 
			|\varPhi^{\epsilon}_t- \bar{\varPhi}_t| 
			\mathrm{d}t\mathrm{d}s\qquad\qquad\quad\quad\\
			&=(1+C)K^2 \int_0^\tau (\tau-t) 
			|\varPhi^{\epsilon}_t- \bar{\varPhi}_t| 
			\mathrm{d}t
			\leq(1+C)K^2 \int_0^\tau 
			|\varPhi^{\epsilon}_t- \bar{\varPhi}_t| 
			\mathrm{d}t.
		\end{align*}
		
		For the term $J_{4}$, by \eqref{bb2}, \eqref{bb'1} and Lemma \ref{Zbbounded}, we obtain
		\begin{align*}
			|J_{4}|
			&\leq \int_0^\tau\int_0^s 
			K(1+|\bar{\varPhi}^{\epsilon}_t|)K|\varPhi^{\epsilon}_t- \bar{\varPhi}_t| 
			\mathrm{d}t\mathrm{d}s
			\leq (1+\bar{C})K^2 \int_0^\tau\int_0^s 
			|\varPhi^{\epsilon}_t- \bar{\varPhi}_t| 
			\mathrm{d}t\mathrm{d}s\qquad\qquad\qquad\\
			&=(1+\bar{C})K^2 \int_0^\tau (\tau-t) 
			|\varPhi^{\epsilon}_t- \bar{\varPhi}_t| 
			\mathrm{d}t
			\leq(1+\bar{C})K^2 \int_0^\tau 
			|\varPhi^{\epsilon}_t- \bar{\varPhi}_t| 
			\mathrm{d}t.
		\end{align*}
		
		For the term $J_{5}$, by exchanging the order of integration and the Lipschitz condition \eqref{bb'1}, we have
		\begin{align*}
			|J_{5}|
			\leq \frac{1}{2} \int_0^\tau 
			|b''(t/\epsilon,\varPhi^{\epsilon}_t)-b''(t/\epsilon,\bar{\varPhi}_t)|\mathrm{d}t
			\leq  \frac{  K}{2}\int_0^\tau |\varPhi^{\epsilon}_t-\bar{\varPhi}_t|\mathrm{d}t.\qquad\qquad\qquad\qquad\qquad\qquad
		\end{align*}
		
		For the term $J_{6}$, by exchanging the order of integration, change of variables and the averaged condition \eqref{acb'}, we obtain
		\begin{align*}
			|J_{6}|&\leq \frac{1}{2} \int_0^\tau\int_0^s 
			|b''(t/\epsilon,\bar{\varPhi}_t)-\bar{b}''(\bar{\varPhi}_t)|\mathrm{d}t\mathrm{d}s
			\leq \frac{1}{2} \int_0^\tau 
			|b''(t/\epsilon,\bar{\varPhi}_t)-\bar{b}''(\bar{\varPhi}_t)|\mathrm{d}t\qquad\qquad\qquad\quad\\
			&\leq \frac{  \epsilon}{2\tau} \int_0^{\frac{\tau}{\epsilon}} 
			|b''(u,\bar{\varPhi}_{u\epsilon})-\bar{b}''(\bar{\varPhi}_{u\epsilon})|\mathrm{d}u
			\leq \frac{  (1+\bar{C})}{2}\phi_2\left(\frac{\tau}{\epsilon}\right).
		\end{align*}
		
		Next, we deal with the term $J_{7}$, by the chain rule $\mathrm{d}b(x,y(x))=\partial_1 b(x,y)+\partial_2 b(x,y)\frac{\mathrm{d}y}{\mathrm{d}x}$, we know that
		\begin{align*}
			J_{7}&
			=\frac{1}{\epsilon}  \int_0^\tau 
			(\tau-t)\partial_tb(t/\epsilon,\varPhi^{\epsilon}_t)\mathrm{d}t
			=\frac{1}{\epsilon}  \int_0^\tau 	(\tau-t)
			\left[\frac{\epsilon\mathrm{d}b(t/\epsilon,\varPhi^{\epsilon}_t)}{\mathrm{d}t}
			-\epsilon b'(t/\epsilon,\varPhi^{\epsilon}_t)\dot{\varPhi}^{\epsilon}_t\right]
			\mathrm{d}t\qquad\quad\\
			&=\int_0^\tau (\tau-t)\mathrm{d}b(t/\epsilon,\varPhi^{\epsilon}_t)\mathrm{d}t
			-\int_0^\tau (\tau-t)b'(t/\epsilon,\varPhi^{\epsilon}_t)\dot{\varPhi}^{\epsilon}_t\mathrm{d}t
			:=J_{71}-J_{72}.
		\end{align*}
		Using integration by parts and changing variables, the term $J_{71}$ becomes
		\begin{align*}
			J_{71}&=\int_0^\tau (\tau-t)\mathrm{d}b(t/\epsilon,\varPhi^{\epsilon}_t)\mathrm{d}t
			=(\tau-t)b(t/\epsilon,\varPhi^{\epsilon}_t)\mid^{t=\tau}_{t=0}
			+\int_0^\tau b(t/\epsilon,\varPhi^{\epsilon}_t)\mathrm{d}t\qquad\qquad\qquad\;\quad\\
			&= -\tau b(0,\varPhi^{\epsilon}_0)
			+\epsilon\int_0^{\frac{\tau}{\epsilon}} b(u,\varPhi^{\epsilon}_{u\epsilon})\mathrm{d}u
			=\tau \frac{\epsilon}{\tau}\int_0^{\frac{\tau}{\epsilon}} b(u,\varPhi^{\epsilon}_{u\epsilon})\mathrm{d}u,
		\end{align*}
		where we assume the standard initial layer condition $b(0,\varPhi^{\epsilon}_0) = 0$. Furthermore, under the assumption that the time-average of the drift vanishes, i.e.,
		$$
		\lim_{T_1\rightarrow \infty}\frac{1}{T_1}\int_{0}^{T_1}b(t,y)\mathrm{d}t=\bar{b}(y)=0,
		$$
		we have $\lim_{\epsilon \rightarrow 0} \sup_{\tau \in [0,1]} |J_{71}| = 0$ as $\epsilon \rightarrow 0$.
		
		For the term $J_{72}$, we utilize the anti-derivative of the spatial derivative with respect to the fast time variable. 
		By the chain rule, we can express the highly oscillatory term as
		$$
		\partial_x b(t/\epsilon, \varPhi^\epsilon_t) = \epsilon \frac{\mathrm{d}}{\mathrm{d}t} \left( \int_0^{t/\epsilon} \partial_x b(r, \varPhi^\epsilon_t) \mathrm{d}r \right) - \epsilon \left( \int_0^{t/\epsilon} \partial_{x}^2 b(r, \varPhi^\epsilon_t) \mathrm{d}r \right) \dot{\varPhi}^\epsilon_t.
		$$
		Substituting this identity directly into $J_{72}$ yields
		\begin{align*}
			J_{72} &= \int_0^\tau (\tau-t) \left\{ \epsilon \frac{\mathrm{d}}{\mathrm{d}t} \left( \int_0^{t/\epsilon} \partial_x b(r, \varPhi^\epsilon_t) \mathrm{d}r \right) - \epsilon \left( \int_0^{t/\epsilon} \partial_{x}^2 b(r, \varPhi^\epsilon_t) \mathrm{d}r \right) \dot{\varPhi}^\epsilon_t \right\} \mathrm{d}t.
		\end{align*}
		Integrating the first term by parts transfers the time derivative away from the highly oscillating anti-derivative
		\begin{align*}
			\epsilon \int_0^\tau (\tau-t) \frac{\mathrm{d}}{\mathrm{d}t} \left( \int_0^{t/\epsilon} \partial_x b(r, \varPhi^\epsilon_t) \mathrm{d}r \right) \mathrm{d}t 
			&= \epsilon \int_0^\tau \left( \int_0^{t/\epsilon} \partial_x b(r, \varPhi^\epsilon_t) \mathrm{d}r \right) \mathrm{d}t,
		\end{align*}
		Therefore, $J_{72}$ can be bounded by
		\begin{align*}
			|J_{72}| &\leq 
			\epsilon \int_0^\tau \left| \int_0^{t/\epsilon} \partial_x b(r, \varPhi^\epsilon_t) \mathrm{d}r \right| \mathrm{d}t 
			+ \epsilon \int_0^\tau (\tau-t) \left| \int_0^{t/\epsilon} \partial_{x}^2 b(r, \varPhi^\epsilon_t) \mathrm{d}r \right| |\dot{\varPhi}^\epsilon_t| \mathrm{d}t.\quad
		\end{align*}
		By Lemma \ref{Psibounded}, growth conditions \eqref{b2} and \eqref{b'2}, we can obtain that there exists a positive constant $C$ such that 
		$\sup_{0\leq \tau\leq 1}|\dot{\varPhi}^{\epsilon}_\tau|\leq C$. 
		Under the averaged assumption $\bar{b}(x) = 0$, its spatial derivatives also have zero mean, ensuring that the integrals $\int_0^\theta \partial_x b(r, x) \mathrm{d}r$ and $\int_0^\theta \partial_{x}^2 b(r, x) \mathrm{d}r$ are bounded uniformly for all $\theta \geq 0$. 
		Therefore, every integral on the right-hand side is finite and the explicit pre-factor $\epsilon$ strictly guarantees that
		$$
		\lim_{\epsilon \rightarrow 0} \sup_{\tau \in [0,1]} |J_{72}| = 0.
		$$
		Taking all these estimations on $J_1$-$J_7$ into consideration and using Gronwall inequatity, we get 
		$$
		\lim_{\epsilon \rightarrow 0}\sup_{\tau \in [0,1]}
		\left|\varPhi^\epsilon_\tau - \bar{\varPhi}_\tau\right|=0.
		$$
		The proof is complete.
	\end{proof}
	
	\begin{remark}
		The conditions imposed on the drift term $b(t,y)$, namely that $b(0,x_0)=0$ and its time average vanishes in the long-time limit for all $y \in \mathbb{R}^d$, are not overly restrictive. 
		Physically, they describe systems driven by transient or zero-mean oscillatory forces without a permanent macroscopic bias. 
		Typical functions satisfying these requirements include
		\begin{enumerate}
			\item Transient decaying drift: $b(t,y) = t e^{-\alpha t} g(y)$ for $\alpha > 0$. The factor $t$ ensures $b(0,x_0)=0$, while the exponential decay guarantees the long-time average converges to zero.
			\item Zero-mean oscillatory drift: $b(t,y) = \sin(\omega t) g(y)$ for $\omega \neq 0$. It trivially satisfies $b(0,x_0)=0$, and the periodic oscillation yields a zero  average over an infinite time horizon.
		\end{enumerate}
	\end{remark}
	\section{$\Gamma$-Convergence of Onsager-Machlup action functionals}\label{GC}
	The present section establishes a formal result regarding the convergence of both the minima and the minimizers within the framework of $\Gamma$-convergence. 
	Specifically, for the SDE \eqref{SDE11}, we investigate the asymptotic behavior, as $\epsilon \to 0$, of the minima and minimizers for the family of functionals $\{F^\epsilon\}_{\epsilon > 0}: H^1(0,1) \to \mathbb{R}$ defined by \begin{equation}\label{Feps}
		F^\epsilon(\varphi^\epsilon, \dot{\varphi}^\epsilon)
		= \int_{0}^{1} \left[ \left(b(t/\epsilon, \varphi^{\epsilon}_t) - \dot{\varphi}^\epsilon_t\right)^2 + b'(t/\epsilon, \varphi^{\epsilon}_t) \right] \mathrm{d}t,
	\end{equation}
	and demonstrate that the sequence $\{F^\epsilon\}_{\epsilon>0}$ $\Gamma$-converges to the homogenized limit functional $F^0$ defined by 
	\begin{equation}\label{F0}
		F^0(\bar{\varphi}, \dot{\bar{\varphi}}) = \int_{0}^{1} \left[ \left(\bar{b}(\bar{\varphi}_t) - \dot{\bar{\varphi}}_t\right)^2 + \bar{b}'(\bar{\varphi}_t) \right] \mathrm{d}t.
	\end{equation}
	The main result of this section is as follows.	
	\begin{theorem}[Convergence of Minima and Minimizers]
		\label{thm:convergence_minima}
		Let $\{F^\epsilon\}_{\epsilon > 0}$ be the family of equi-mildly coercive functionals \eqref{Feps} and $F^0$ defined by \eqref{F0} be its $\Gamma$-limit in the weak topology of $H^1(0,1)$ and consequently in the strong topology of $L^2(0,1)$. Under the assumptions $\textbf{(A1)-(A2)}$, the following assertions hold:
		\begin{enumerate}
			\item \textbf{Convergence of Minima:} The sequence of the values of the minima converges to the minimum of the limit functional, i.e.,
			\begin{equation*}
				\lim_{\epsilon \to 0} \left( \inf_{\varphi \in H^1} F^\epsilon(\varphi) \right) = \min_{\varphi \in H^1} F^0(\varphi).
			\end{equation*}
			
			\item \textbf{Convergence of Minimizers:} Let $\{\varphi^\epsilon\}_{\epsilon > 0}$ be a sequence of almost-minimizers for $F^\epsilon$, satisfying $\lim_{\epsilon \to 0} (F^\epsilon(\varphi^\epsilon) - \inf F^\epsilon) = 0$. Every cluster point $\bar{\varphi}$ of the sequence $\{\varphi^\epsilon\}$ with respect to the weak topology of $H^1(0,1)$ is a global minimizer of the homogenized functional $F^0$.
			\item \textbf{Uniqueness:} If $F^0$ possesses a unique global minimizer $\bar{\varphi}$, then the entire sequence $\{\varphi^\epsilon\}$ weakly converges to $\bar{\varphi}$ in $H^1(0,1)$ as $\epsilon \to 0$.
		\end{enumerate}
	\end{theorem}
	The theorem follows directly from the fundamental result of $\Gamma$-convergence (Lemma \ref{min_convergence}) by combining the $\Gamma$-convergence and the equi-coercivity,
	which we will establish in Lemma \ref{gamma4.1} and Lemma \ref{equi-mildly coercive}  separately. 
	Consequently, our first goal is to establish the $\Gamma$-convergence of the sequence.
	\begin{lemma}\label{gamma4.1} 
		Let $F^\epsilon$ be defined as above. 
		Under the assumptions $\textbf{(A1)-(A2)}$,  
		the sequence of functionals $F^\epsilon$ $\Gamma$-converges in the weak topology of $H^1(0,1)$ (and consequently, the strong topology of $L^2(0,1)$) to the functional $F^0$.
	\end{lemma}  
	\begin{proof}
		By Definition \ref{gamma}, the proof consists of verifying the two standard conditions of $\Gamma$-convergence, that is, the liminf inequality for any weakly convergent sequence, and the limsup inequality asserting the existence of a recovery sequence.
		
		\textit{Step 1: The Liminf Inequality.}
		Let $\{\varphi^\epsilon\}$ be an arbitrary sequence such that $\varphi^\epsilon \rightharpoonup \bar{\varphi}$ weakly in $H^1(0,1)$. 
		By the Rellich-Kondrachov compactness theorem(see Theorem 1 on page 272 of \cite{evans2022partial}), $\varphi^\epsilon \to \bar{\varphi}$ strongly in $L^2(0,1)$ and uniformly in $C([0,1])$. 
		Consequently, the sequence is uniformly bounded, meaning there exists a constant $M > 0$ such that $\|\varphi^\epsilon\| \leq M$ for all $\epsilon > 0$.
		
		We rewrite the main quadratic term of the functional by adding and subtracting $\bar{b}(\varphi^\epsilon_t)$
		\begin{align}\label{quadratic}
			\left(b(t/\epsilon, \varphi^{\epsilon}_t) - \dot{\varphi}^\epsilon_t\right)^2 
			&= \left( (\bar{b}(\varphi^\epsilon_t) - \dot{\varphi}^\epsilon_t) + (b(t/\epsilon, \varphi^\epsilon_t) - \bar{b}(\varphi^\epsilon_t)) \right)^2 \nonumber\\
			&= \left( \bar{b}(\varphi^\epsilon_t) - \dot{\varphi}^\epsilon_t \right)^2 + \left( b(t/\epsilon, \varphi^\epsilon_t) - \bar{b}(\varphi^\epsilon_t) \right)^2 \nonumber\\
			&\quad + 2 \left( \bar{b}(\varphi^\epsilon_t) - \dot{\varphi}^\epsilon_t \right)\left( b(t/\epsilon, \varphi^\epsilon_t) - \bar{b}(\varphi^\epsilon_t) \right).
		\end{align}
		
		Let us analyze the discrepancy term $D^\epsilon(t) := b(t/\epsilon, \varphi^\epsilon_t) - \bar{b}(\varphi^\epsilon_t)$. Utilizing the strong averaging condition \eqref{ac} and the uniform $L^\infty$ bound $M$, we estimate its $L^1$ norm	
		\begin{equation*}
			\int_0^1 |D^\epsilon(t)| \mathrm{d}t = \epsilon \int_0^{1/\epsilon} |b(s, \varphi^\epsilon_{\epsilon s}) - \bar{b}(\varphi^\epsilon_{\epsilon s})| \mathrm{d}s \leq \epsilon \left( \frac{1}{\epsilon} \phi_1(1/\epsilon) (1 + M) \right) = \phi_1(1/\epsilon)(1+M).
		\end{equation*}
		As $\epsilon \to 0$, $1/\epsilon \to \infty$, which yields $\phi_1(1/\epsilon) \to 0$. 
		Hence, $D^\epsilon \to 0$ strongly in $L^1(0,1)$.
		Given the $L^\infty$ bound on $\varphi^\epsilon$,
		$D^\epsilon$ is uniformly bounded in $L^\infty$, implying that $D^\epsilon \to 0$ strongly in $L^2(0,1)$ as well.
		Since $\varphi^\epsilon \rightharpoonup \bar{\varphi}$ weakly in $H^1$, $\dot{\varphi}^\epsilon$ is bounded in $L^2$. Thus, under the grown condition \eqref{bb2} of the function $\bar{b}$, we can obtain that $\dot{\varphi}^\epsilon_t - \bar{b}(\varphi^\epsilon_t)$ is bounded in $L^2$ by some constant $M_1$. 
		
		Now, substituting the expansion into $F^{\epsilon}$, by the Cauchy-Schwarz inequality, we have
		\begin{align*}\label{eq:expansion}
			&\quad F^{\epsilon}(\varphi^\epsilon) 
			= \int_0^1 |\dot{\varphi}^\epsilon_t - \bar{b}(\varphi^\epsilon_t)|^2 \mathrm{d}t 
			+ \int_0^1 |D^\epsilon(t)|^2 \mathrm{d}t 
			- 2 \int_0^1 \left( 
			\dot{\varphi}^\epsilon_t - \bar{b}(\varphi^\epsilon_t) 
			\right) D^\epsilon(t) \mathrm{d}t 
			+ \int_0^1 b'(t/\epsilon, \varphi^{\epsilon}_t) \mathrm{d}t\nonumber\\
			&\geq \int_0^1\! |\dot{\varphi}^\epsilon_t - \bar{b}(\varphi^\epsilon_t)|^2 \mathrm{d}t
			\!+\! \int_0^1 \!|D^\epsilon(t)|^2 \mathrm{d}t 
			\!-\!2\left(\!\int_0^1\!|\dot{\varphi}^\epsilon_t - \bar{b}(\varphi^\epsilon_t)|^2\mathrm{d}t \!\right)^{\frac{1}{2}}
			\left(\!\int_0^1\! |D^\epsilon(t)|^2\mathrm{d}t \!\right)^{\frac{1}{2}}
			\!+\!\int_0^1\! b'(t/\epsilon, \varphi^{\epsilon}_t) \mathrm{d}t
			\nonumber\\
			&\geq \int_0^1 |\dot{\varphi}^\epsilon_t - \bar{b}(\varphi^\epsilon_t)|^2 \mathrm{d}t
			+ \int_0^1 |D^\epsilon(t)|^2 \mathrm{d}t 
			-2\sqrt{M_1}\left(\int_0^1 |D^\epsilon(t)|^2\mathrm{d}t \right)^{\frac{1}{2}}
			+ \int_0^1 b'(t/\epsilon, \varphi^{\epsilon}_t) \mathrm{d}t.
		\end{align*}
		
		Take the limit inferior of the above inequality, the second and third term tend to zero as $\epsilon\rightarrow 0$. 
		By the uniform convergence of $\varphi^\epsilon \to \bar{\varphi}$ and the dominated convergence theorem, we have
		$$
		\lim_{\epsilon \to 0} \int_0^1 b'(t/\epsilon, \varphi^{\epsilon}_t) \mathrm{d}t = \int_0^1 \bar{b}'(\bar{\varphi}_t) \mathrm{d}t,
		$$
		which yeilds 
		\begin{equation*}
			\liminf_{\epsilon \to 0} F^{\epsilon}(\varphi^\epsilon) 
			\geq \liminf_{\epsilon \to 0} 
			\int_0^1 \left( \dot{\varphi}^\epsilon_t - \bar{b}(\varphi^\epsilon_t) \right)^2 \mathrm{d}t 
			+ \int_0^1 \bar{b}'(\bar{\varphi}_t) \mathrm{d}t.
		\end{equation*}
		Since $\varphi^\epsilon \to \bar{\varphi}$ strongly in $L^2$ and $\dot{\varphi}^\epsilon \rightharpoonup \dot{\bar{\varphi}}$ weakly in $L^2$, the sequence $\dot{\varphi}^\epsilon_t - \bar{b}(\varphi^\epsilon_t)$ converges weakly in $L^2$ to $\dot{\bar{\varphi}}_t - \bar{b}(\bar{\varphi}_t)$. By the weak lower semicontinuity of the $L^2$ norm,
		\begin{equation*}
			\liminf_{\epsilon \to 0} \int_0^1 \left( \dot{\varphi}^\epsilon_t - \bar{b}(\varphi^\epsilon_t) \right)^2 \mathrm{d}t \geq \int_0^1 \left( \dot{\bar{\varphi}}_t - \bar{b}(\bar{\varphi}_t) \right)^2 \mathrm{d}t.
		\end{equation*}
		Therefore, $\liminf_{\epsilon \to 0} F^{\epsilon}(\varphi^\epsilon) \geq F^0(\bar{\varphi})$, which establishes the liminf inequality.
		
		\textit{Step 2: The Limsup Inequality (Recovery Sequence).} 
		To establish the limsup inequality for any given $\bar{\varphi} \in H^1$, we must construct a recovery sequence $\varphi^\epsilon \to \bar{\varphi}$ such that $\limsup_{\epsilon \to 0} F^{\epsilon}(\varphi^\epsilon) \leq F^0(\bar{\varphi})$. 
		Because the condition \eqref{ac} dictates a strong convergence of the oscillatory coefficients to their mean, we do not need to construct a sequence with microscopic oscillations. We can simply choose the constant sequence
		\begin{equation*}
			\varphi^\epsilon(t) \equiv \bar{\varphi}(t) \quad \text{for all } \epsilon > 0.
		\end{equation*}
		Substituting this sequence into the functional gives
		\begin{equation*}
			F^{\epsilon}(\varphi^\epsilon)=F^{\epsilon}(\bar{\varphi}) 
			= \int_{0}^{1} \left[ \left(\dot{\bar{\varphi}}_t - b(t/\epsilon, \bar{\varphi}_t)\right)^2 + b'(t/\epsilon, \bar{\varphi}_t) \right]\mathrm{d}t.
		\end{equation*}
		We expand the quadratic term exactly as \eqref{quadratic} in Step 1. Since $D^\epsilon(t) = b(t/\epsilon, \bar{\varphi}_t) - \bar{b}(\bar{\varphi}_t) \to 0$ strongly in $L^2(0,1)$, the term $\int_0^1 |D^\epsilon|^2 \mathrm{d}t$ and the cross-term both converge to zero as $\epsilon\rightarrow 0$. 
		Therefore, we have
		\begin{align*}
			\limsup_{\epsilon \to 0} F^{\epsilon}(\bar{\varphi}) &\leq \lim_{\epsilon \to 0} \int_{0}^{1} \left[ \left(\dot{\bar{\varphi}}_t - b(t/\epsilon, \bar{\varphi}_t)\right)^2 + b'(t/\epsilon, \bar{\varphi}_t) \right]\mathrm{d}t \\
			&\leq \int_{0}^{1} \left[ \left(\dot{\bar{\varphi}}_t - \bar{b}(\bar{\varphi}_t)\right)^2 + \bar{b}'(\bar{\varphi}_t) \right]\mathrm{d}t = F^0(\bar{\varphi}).
		\end{align*}
		This proves the limsup inequality. 
		
		Combining this with the liminf inequality, we conclude that $\Gamma\text{-}\lim_{\epsilon \to 0} F^{\epsilon} = F^0$ in the weak topology of $H^1$, and consequently, by the Rellich-Kondrachov compactness theorem, it also holds in the strong topology of $L^2(0,1)$.
		The proof is completed.
	\end{proof}
	
	Now our second goal is to prove the equi-coerciveness. 
	\begin{lemma}\label{equi-mildly coercive}
		Under the assumptions $\textbf{(A1)-(A2)}$, the family of functionals $\{F^\epsilon\}_{\epsilon > 0}$ is equi-coercive on $H^1(0,1)$. That is, for any $L \in \mathbb{R}$, the common sublevel set
		\begin{equation*}
			\mathcal{K}_L := \bigcup_{\epsilon > 0} \left\{ \varphi^\epsilon \in X \;\middle|\; F^\epsilon(\varphi^\epsilon) \leq L\right\}
		\end{equation*}
		is precompact in $H^1(0,1)$. Consequently, $\{F^\epsilon\}_{\epsilon > 0}$ is equi-coercive.
	\end{lemma}
	
	\begin{proof} 
		For any fixed $L \in \mathbb{R}$. Let $\{\varphi^\epsilon\}_{\epsilon > 0}$ be a family of absolutely continuous paths in $\mathcal{K}_L$. 
		From the proof of Lemma \ref{gamma4.1}, there exists a constant $M > 0$ such that $\|\varphi^\epsilon\|_{L^\infty} \leq M$ for all $\epsilon > 0$, so the family $\mathcal{K}_L$ is uniformly bounded.
		
		On the other hand, by the growth condition \eqref{b'2}, the divergence term $b'(t,y)$ satisfies $\left| b'(t,y) \right| \leq K(1 + \left| y \right|)$, which means
		\begin{equation*}
			b'(t/\epsilon, \varphi^\epsilon_t) \ge - \left| b'(t/\epsilon, \varphi^\epsilon_t) \right| \ge - K(1 + M) =: -C_1,
		\end{equation*}
		where $C_1 > 0$ is a constant independent of $\epsilon$ and $t$. Integrating this over $[0,1]$ yields
		\begin{equation*} \label{eq:b_prime_bound}
			\int_{0}^{1} b'(t/\epsilon, \varphi^\epsilon_t) \mathrm{d} t \ge -C_1.
		\end{equation*}	
		Since $F^\epsilon(\varphi^\epsilon) \le L$, we have
		\begin{equation} \label{eq:diff_square_bound}
			\int_{0}^{1} \left| \dot{\varphi}^\epsilon_t - b(t/\epsilon, \varphi^{\epsilon}_t) \right|^2 \mathrm{d}t 
			=F^\epsilon(\varphi^\epsilon)-\int_{0}^{1} b'(t/\epsilon, \varphi^\epsilon_t) \mathrm{d}t 
			\leq L + C_1.
		\end{equation}
		Using the elementary algebraic inequality $\left| u \right|^2 \le 2\left| u-v \right|^2 + 2\left| v \right|^2$, the growth condition \eqref{b2}, \eqref{eq:diff_square_bound} and the fact that $\left\lVert \varphi^\epsilon \right\rVert_{L^\infty} \le M$,
		we can obtain
		\begin{equation*} \label{eq:kinetic_energy}
			\int_{0}^{1} \left| \dot{\varphi}^\epsilon_t \right|^2 \mathrm{d} t 
			\leq 2 \int_{0}^{1} \left| \dot{\varphi}^\epsilon_t - b(t/\epsilon, \varphi^{\epsilon}_t) \right|^2 \mathrm{d} t 
			+ 2 \int_{0}^{1} \left| b(t/\epsilon, \varphi^\epsilon_t) \right|^2 \mathrm{d} t
			\leq 2(L + C_1) + 2K^2(1+M)^2.
		\end{equation*}
		Thus, we conclude that the family of paths $\{\varphi^\epsilon\}_{\epsilon > 0}$ is uniformly bounded.
		
		By virtue of the Banach-Alaoglu theorem, in a reflexive space, norm boundedness directly implies weak precompactness. Furthermore, we appeal to the Eberlein-\v{S}mulian theorem, which guarantees that for Banach spaces, weak precompactness is strictly equivalent to sequential weak precompactness. 
		Consequently, any sequence extracted from the sublevel set $\mathcal{K}_L$ necessarily admits a subsequence that converges weakly in $H^1(0,1)$. This  demonstrates that $\mathcal{K}_L$ is weakly precompact in $H^1(0,1)$, which completes the proof.		
	\end{proof}
	\begin{remark}
		Since the underlying space $X = H^1(0,1)$ is an infinite-dimensional Hilbert space, the compactness required in Definition \ref{coerciveness} is understood in the sense of the weak topology. 
		As we have demonstrated that the sublevel set $\mathcal{K}_L$ is weakly precompact, we can define $K$ as the weak closure of $\mathcal{K}_L$. Thus, $K$ serves as the common weakly compact set, which rigorous confirms that the family $\{F_\epsilon\}_{\epsilon > 0}$ is equi-mildly coercive with respect to the weak topology of $H^1(0,1)$.
	\end{remark}
	\section{Numerical experiments}\label{NE}
	In this section, we provide numerical examples to illustrate the theoretical results established in Theorem \ref{main} and Theorem \ref{mainn}. We consider a one-dimensional stochastic system governed by a double-well potential function $V(x) = -\frac{1}{2}x^2 + \frac{1}{4}x^4$. The potential wells are located at $x = \pm 1$, representing the metastable states, separated by a potential barrier at $x = 0$. We investigate the most probable transition paths between these metastable states under different scaling regimes.
	
	$\bullet$ Case 1: Proximity between the most probable transition paths $Z^{\epsilon}$ and $\bar{Z}^{\epsilon}$.
	
	First, we consider the time-inhomogeneous SDE:
	\begin{equation*}
		\mathrm{d}X^{\epsilon}_t = \epsilon \sin (2\pi t) [X^{\epsilon}_t - (X^{\epsilon}_t)^3] \mathrm{d}t + \sqrt{\epsilon} \mathrm{d}B_t, \quad X^{\epsilon}_0 = 1.
	\end{equation*}
	The corresponding derivatives of the drift are given by $b'(t,x) = \sin(2\pi t)(1-3x^2)$, $b''(t,x) = -6x\sin(2\pi t)$, and $\partial_t b(t,x) = 2\pi(x-x^3)\cos (2\pi t)$. 
	By OM theory, we konow that the MPTP $Z^\epsilon_t$ is the solution to the following Euler--Lagrange equation:
	\begin{equation*}
		\begin{aligned}
			\ddot{Z}^\epsilon_t &= \epsilon^2 b(t,Z^{\epsilon}_t)b'(t,Z^{\epsilon}_t) + \frac{\epsilon^2}{2} b''(t,Z^{\epsilon}_t) + \epsilon \partial_t b(t,Z^{\epsilon}_t) \\
			&= \epsilon^2 \sin^2(2\pi t) [Z^{\epsilon}_t - (Z^{\epsilon}_t)^3][1-3(Z^{\epsilon}_t)^2] - 3\epsilon^2 Z^{\epsilon}_t \sin(2\pi t) + 2\pi\epsilon [Z^{\epsilon}_t - (Z^{\epsilon}_t)^3] \cos(2\pi t),
		\end{aligned}
	\end{equation*}
	subject to the boundary conditions $Z^\epsilon_0 = 1$ and $Z^\epsilon_1 = -1$.
	\begin{figure}[htbp]
		\centering
		\includegraphics[width=1\textwidth]{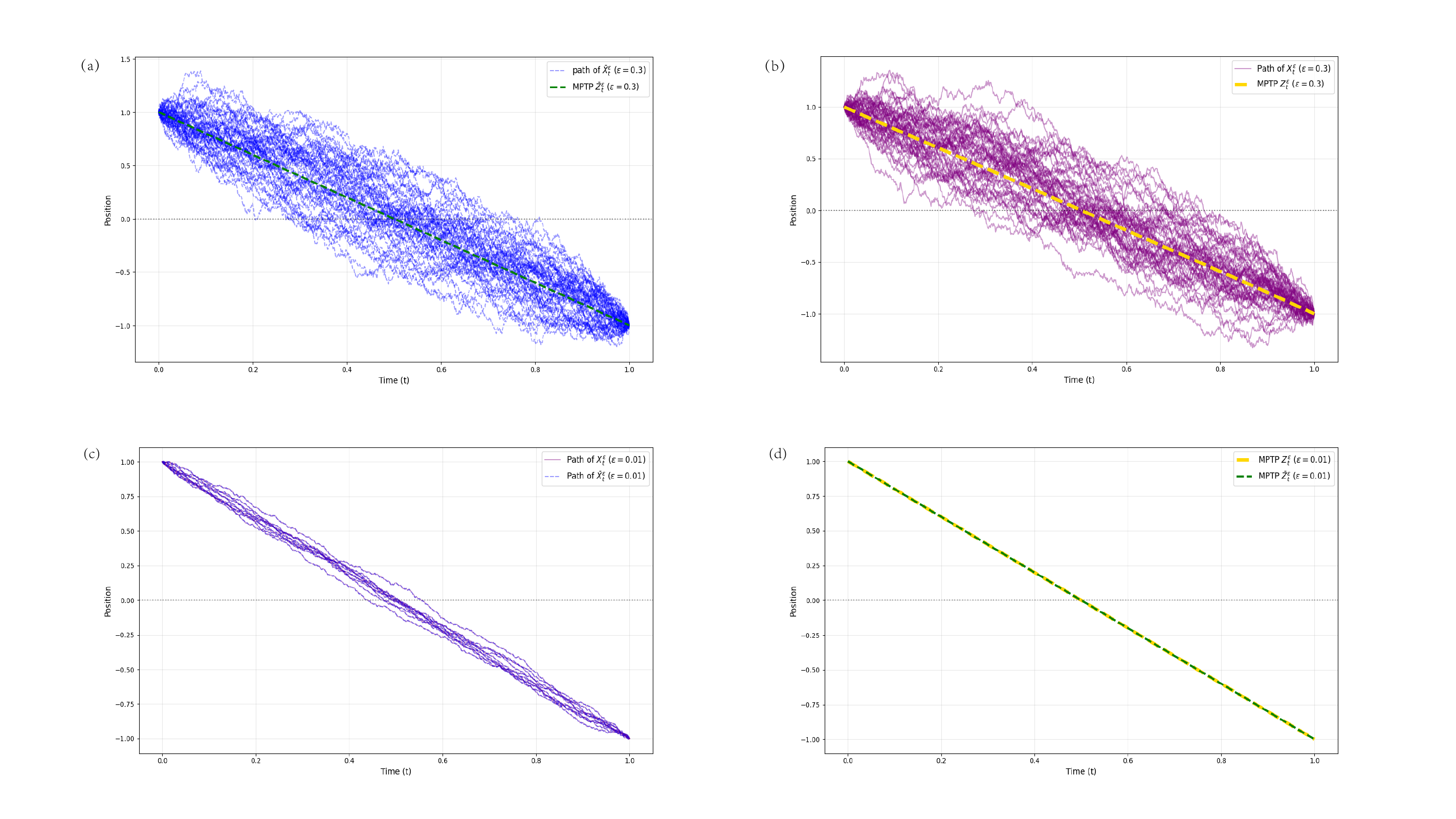}
		\caption{Comparisons of SDEs ($X^{\epsilon}, \bar{X}^{\epsilon}$) and the most probable transition paths ($Z^{\epsilon}, \bar{Z}^{\epsilon}$).
			(a) and (b) present the sample paths of the SDEs (blue and purple lines) and their respective MPTPs (green and yellow lines) for both the original and proximate averaged systems with $\epsilon=0.3$;
			(c) illustrates the dynamics of the original solution $X^{\epsilon}$ (blue lines)  against the proximate averaged solution $\bar{X}^{\epsilon}$ (purple lines) with $\epsilon=0.01$ ; 
			(d) provides a direct comparison between the two MPTPs, $Z^{\epsilon}$  (green line) and $\bar{Z}^{\epsilon}$ (yellow line) with $\epsilon=0.01$.}
		\label{sde1}
	\end{figure}
	On the other hand, for small $\epsilon$, the proximate averaged dynamics can be described by a time-homogeneous SDE:
	\begin{equation*}
		\mathrm{d}\bar{X}^{\epsilon}_t = \bar{b}(\bar{X}^{\epsilon}_t)\mathrm{d}t + \sqrt{\epsilon} \mathrm{d}B_t = \sqrt{\epsilon} \mathrm{d}B_t, \quad \bar{X}^{\epsilon}_0 = 1,
	\end{equation*}
	where the averaged drift $\bar{b}(x) = \int_0^1 \sin(2\pi t)(x-x^3) \mathrm{d}t = 0$.
	Also, the MPTP $\bar{Z}^{\epsilon}_t$ satisfies $\ddot{\bar{Z}}^{\epsilon}_t = 0$ with the same boundary conditions. 
	Since all the assumptions $\textbf{(A1)-(A4)}$ are satisfied, Theorem \ref{main} holds, that is, there exists a constant $C > 0$ such that
	\begin{equation*}
		\sup_{t \in [0, 1]} |Z^\epsilon_t - \bar{Z}^\epsilon_t| \leq C \epsilon.
	\end{equation*}
	
	Now we conduct numerical simulations to investigate the trajectories of the SDEs ($X^{\epsilon}, \bar{X}^{\epsilon}$) and their corresponding MPTPs ($Z^{\epsilon}, \bar{Z}^{\epsilon}$) 
	with the initial state $1$ and terminal state $-1$, the results are summarized in Figure \ref{sde1}. Specifically,  
	Figure \ref{sde1}(a) and (b) present the comparison between the sample paths of the SDEs and their respective MPTPs for both the original and proximate averaged systems with $\epsilon=0.3$;
	Figure \ref{sde1}(c) illustrates the dynamics of the original solution $X^{\epsilon}$ against the proximate averaged solution $\bar{X}^{\epsilon}$ with $\epsilon=0.01$, 
	while Figure \ref{sde1}(d) provides a direct comparison between the two MPTPs, $Z^{\epsilon}$ and $\bar{Z}^{\epsilon}$ with $\epsilon=0.01$. 
	The simulation results exhibit a high degree of consistency between the original MPTP and the proximate averaged one.
	
	$\bullet$ Case 2: Convergence of the most probable transition paths $\varPhi^{\epsilon}$.
	\begin{figure}[htbp]
		\centering
		\includegraphics[width=1\textwidth]{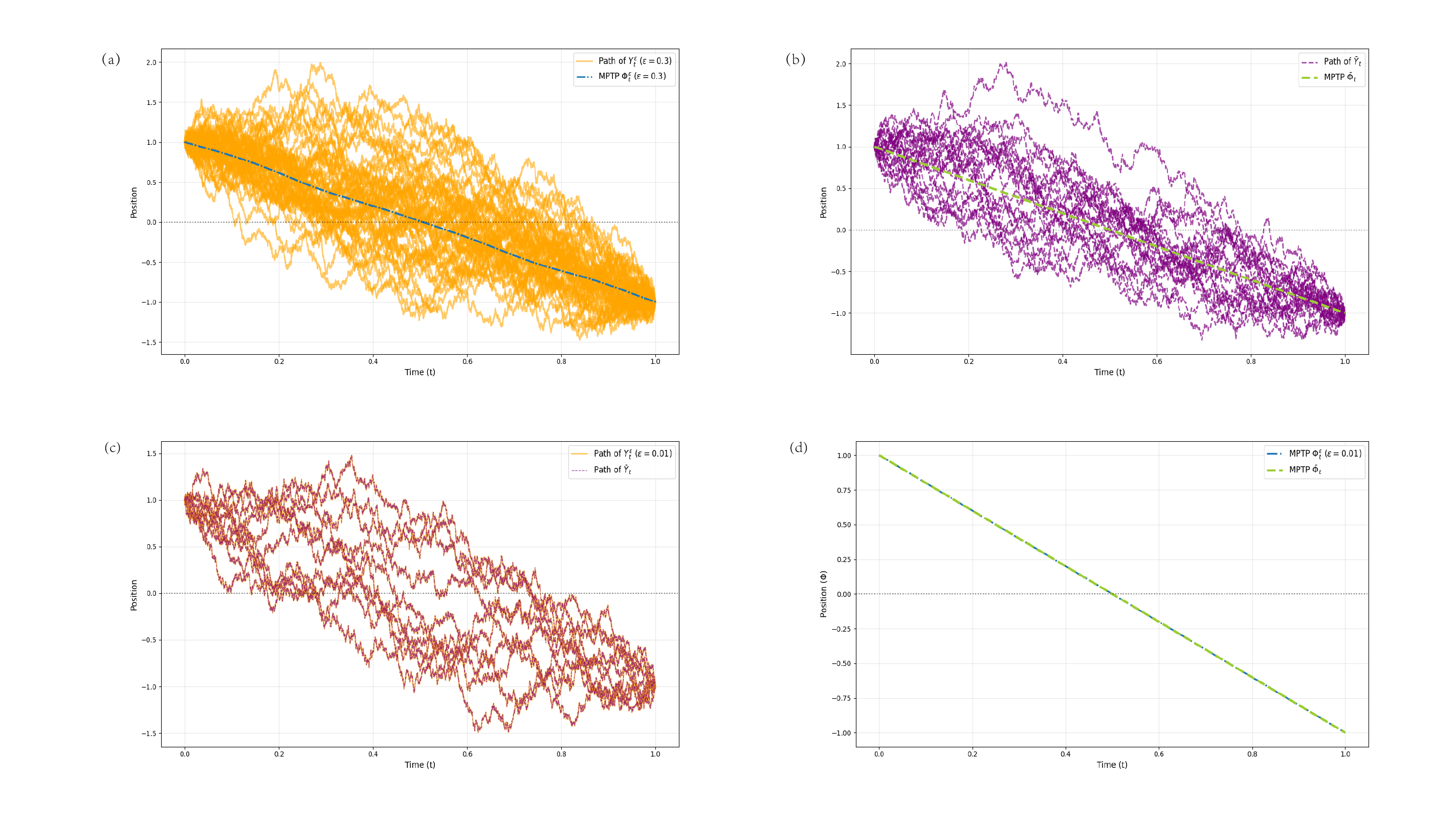}
		\caption{Comparisons of SDEs ($Y^{\epsilon}, \bar{Y}$) and the most probable transition paths ($\varPhi^{\epsilon}, \bar{\varPhi}$). 
			(a) depicts the stochastic sample paths of the original system $Y^{\epsilon}$ (orange lines) and the MPTP $\varPhi^{\epsilon}$ (blue line) with $\epsilon=0.3$;
			(b) depicts the stochastic sample paths of the averaged system $\bar{Y}$ (purple lines) alongside the MPTPs $\bar{\varPhi}$ (yellowgreen line);
			(c) provides a parallel comparison between the original trajectory $Y^{\epsilon}$ (orange lines) with $\epsilon=0.01$ and its averaged counterpart $\bar{Y}$ (purple lines);
			(d) presents a direct contrast between the two MPTPs (blue and yellowgreen lines) with $\epsilon=0.01$.}
		\label{sde2}
	\end{figure}
	
	Next, we consider the regime where the oscillation frequency increases as $\epsilon$ vanishes, leading to a classical averaging scenario. The SDE is given by 
	\begin{equation*} 
		\mathrm{d}Y^{\epsilon}_t = \sin \left(\frac{2\pi t}{\epsilon}\right) [Y^{\epsilon}_t - (Y^{\epsilon}_t)^3] \mathrm{d}t + \mathrm{d}B_t, \quad Y^{\epsilon}_0 = 1,
	\end{equation*}
	with the terminal constraint $Y^{\epsilon}_1 = -1$. 
	The corresponding homogenized averaged system is
	\begin{equation*}
		\mathrm{d}\bar{Y}_t = \mathrm{d}B_t, \quad \bar{Y}_0 = 1, \quad \bar{Y}_1 = -1.
	\end{equation*}
	
	For a fixed $\epsilon$, the most probable path $\varPhi^{\epsilon}_t$ satisfies the Euler-Lagrange equation
	\begin{equation*}
		\ddot{\varPhi}^\epsilon_t = \sin^2\left(\frac{2\pi t}{\epsilon}\right) [\varPhi^{\epsilon}_t - (\varPhi^{\epsilon}_t)^3] [1 - 3(\varPhi^{\epsilon}_t)^2] - 3\varPhi^{\epsilon}_t \sin\left(\frac{2\pi t}{\epsilon}\right) + \frac{2\pi}{\epsilon} [\varPhi^{\epsilon}_t - (\varPhi^{\epsilon}_t)^3] \cos \left(\frac{2\pi t}{\epsilon}\right),
	\end{equation*}
	with boundary conditions $\varPhi^{\epsilon}_0 = 1, \varPhi^{\epsilon}_1 = -1$. 
	Since the coefficients satisfy the conditions in Thoerem \ref{mainn}, it follows that the limit path $\bar{\varPhi}_t$ converges to the solution of $\ddot{\bar{\varPhi}}_t ={\tiny } 0$, as $\epsilon \to 0$, which is the MPTP of the averaged SDE. 
	
	The numerical results are presented in Figure \ref{sde2}. 
	Specifically, Figure \ref{sde2}(a) depicts the stochastic sample paths of the original system $Y^{\epsilon}$ and the MPTP $\varPhi^{\epsilon}$ with $\epsilon=0.3$; 
	Figure \ref{sde2}(b) depicts the stochastic sample paths of the averaged system $\bar{Y}$ alongside the MPTP $\bar{\varPhi}$;
	the two plots illustrate how the MPTP serves as the deterministic backbone of the stochastic transition.
	Figure \ref{sde2}(c) provides a parallel comparison between the original trajectory $Y^{\epsilon}$ with $\epsilon=0.01$ and its averaged counterpart $\bar{Y}$. 
	Figure \ref{sde2}(d) presents a direct contrast between $\varPhi^{\epsilon}$ with $\epsilon=0.01$ and $\bar{\varPhi}$.
	Consistent with our theoretical predictions 	
	$$\lim\limits_{\epsilon \rightarrow 0}\sup_{t \in [0, 1]}
	|\varPhi^\epsilon_t - \bar{\varPhi}_t|=0,$$
	the results exhibit a remarkable alignment between the original MPTP and the averaged MPTP. This agreement demonstrates that the averaged path provides a robust approximation of the most probable transition mechanism, effectively filtering out the high-frequency fluctuations.
	\section{Discussion}\label{DS}
	In this paper, we established proximity and convergence results for the most probable transition paths of a class of multiscale stochastic systems in the space of continuous functions equipped with the supremum norm (Theorems \ref{main} and \ref{mainn}), 
	using the calculus of variations and the Euler-Lagrange equation.
	We further investigated the corresponding Onsager-Machlup functionals from a variational perspective using the framework of $\Gamma$-convergence. 
	By proving the equi-coerciveness of the functional family, we obtained the convergence of the minima and associated minimizers (Theorem \ref{thm:convergence_minima}). 
	As a consequence, the minimizers converge weakly in the Cameron-Martin space $H^1$ and strongly in $L^2$.
	These results provide a rigorous variational justification for the convergence of the most probable transition paths under averaging and establish a connection between the most probable transition paths, averaging principles and $\Gamma$-convergence theory in the study of stochastic dynamical systems.
	Moreover, this work helps provide a rigorous theoretical foundation for the approximation of the most probable transition paths in multiscale stochastic systems and supports the use of averaged models in the study of noise-induced transitions.
	
	{\bf Acknowledgement}
	This work was partly supported by the NSFC International Collaboration Fund for Creative Research Teams (Grant W2541005), 
	the Guangdong Provincial Key Laboratory of Mathematical and Neural Dynamical Systems (Grant 2024B1212010004), 
	the Cross Disciplinary Research Team on Data Science and Intelligent Medicine (2023KCXTD054), 
	the Guangdong-Dongguan Joint Research Fund (Grant 2023A1515140016), 
	Dongguan Key Laboratory for Data Science and Intelligent Medicine, and Natural Science Foundation of Guangdong 
	(Grant NO. 2025A1515011188).
	
	{\bf Data Availability} 
	Data sharing is not applicable to this article as no datasets were generated or analysed during the current study.
	
		 \section*{Declarations}
	{\bf Conflicts of interest} 
	The authors declare no conflict of interests.
	\bibliographystyle{plain}
	\bibliography{RefMPTP}
	
\end{document}